\documentclass[conf]{new-aiaa-custom}

\usepackage[utf8]{inputenc}

\usepackage{graphicx}
\usepackage{amsmath}
\usepackage{amssymb}
\usepackage{bm}
\usepackage[version=4]{mhchem}
\usepackage{longtable,tabularx}

\usepackage{physics} 
\AtBeginDocument{\RenewCommandCopy\qty\SI}

\usepackage{mathtools} 

\usepackage{svg}

\usepackage{smartdiagram}
\usetikzlibrary{trees}
\usepackage{pgfplots}
\pgfplotsset{compat=1.5}

\usepackage{markdown}

\usepackage{booktabs}

\usepackage{hhline}

\usepackage[font={small,bf},labelfont=bf]{caption}

\usepackage{float}

\usepackage{algorithm}
\usepackage{algpseudocode}
\usepackage{setspace}

\titlespacing\section{0pt}{-0.1in}{0.05in}
\titlespacing\subsection{0pt}{-0.075in}{0pt}
\titlespacing\subsubsection{0pt}{0pt}{0pt}

\usepackage{accents}

\usepackage{enumerate}

\usepackage{tcolorbox}
\tcbuselibrary{skins}

\newtcolorbox{mybox}
{
  enhanced jigsaw,
  colframe=black,
  colback=white,
  drop shadow=black!50!white,
  boxrule=0.75pt
}

\newtcolorbox{mybox2}
{
  enhanced jigsaw,
  colframe=black,
  colback=white,
  drop shadow=black!50!white,
  boxrule=0.75pt,
  hbox
}
\newcommand{\behcet}{Beh\c{c}et~A\c{c}\i kme\c{s}e}

\definecolor{darts}{HTML}{009670}
\definecolor{bricks}{HTML}{E74C3C}
\definecolor{fista}{HTML}{075187}
\definecolor{fistagray}{HTML}{929598}
\definecolor{steelblue}{HTML}{4682b4}
\definecolor{goldenrod}{HTML}{daa520}

\newcommand{\q}[1]{#1}
\newcommand{\dq}[1]{\boldsymbol{#1}}

\newcommand{\I}{\mathcal{I}}
\newcommand{\B}{\mathcal{B}}

\newcommand{\R}{\mathbb{R}}

\newcommand{\D}{\mathbb{D}}

\newcommand{\defeq}{\vcentcolon=}

\newcommand{\inone}{\in \mathbb{R}}

\newcommand{\inthree}{\in \mathbb{R}^{3}}
\newcommand{\infour}{\in \mathbb{R}^{4}}
\newcommand{\insix}{\in \mathbb{R}^{6}}
\newcommand{\ineight}{\in \mathbb{R}^{8}}

\newcommand{\infifteen}{\in \mathbb{R}^{15}}

\newcommand{\seco}{\textsc{s}\begin{footnotesize}e\end{footnotesize}\textsc{co}}
\newcommand{\sclerp}{\textsc{s}\begin{footnotesize}c\end{footnotesize}\textsc{lerp}}

\newcommand{\ptr}{\textsc{ptr}}
\newcommand{\dqg}{\textsc{dqg}}

\newcommand{\mosek}{\textsc{mosek}}
\newcommand{\pipg}{\textsc{pipg}}
\newcommand{\pipgc}{\textsc{pipg}\textsubscript{\tiny custom}}

\newcommand{\bsocp}{\textsc{bsocp}}

\newcommand{\cprs}{\textsc{cprs}}
\newcommand{\ecos}{\textsc{ecos}}

\newcommand{\gurobi}{\textsc{gurobi}}

\newcommand{\ctscvx}{{\scalebox{1.1}{\textsc{{\scalebox{0.73}{ct-}}sc{\scalebox{0.73}{vx}}}}}}
\newcommand{\autoscvx}{{\scalebox{1.1}{\textsc{{\scalebox{0.73}{Auto-}}sc{\scalebox{0.73}{vx}}}}}}

\newcommand{\lcvx}{\textsc{lc}\begin{scriptsize}\textsc{vx}\end{scriptsize}}

\newcommand{\blkdiag}{\operatorname{blkdiag}}

\newcommand{\chol}{\operatorname{chol}}

\newcommand{\goldenrod}[1]{\textcolor{goldenrod}{#1}}
\newcommand{\dx}[1]{\textcolor{darts}{#1}}
\newcommand{\rx}[1]{\textcolor{bricks}{#1}}
\newcommand{\fx}[1]{\textcolor{fista}{#1}}

\newcommand{\range}[2]{#1\!:\!#2} 

\newcommand{\circdot}[1]{\accentset{\circ}{#1}} 

\newcommand{\f}[2]{#1\!\left(#2\right)} 

\let\norm\undefined 
\DeclarePairedDelimiter\norm{\|}{\|}

\let\OLDthebibliography\thebibliography
\renewcommand\thebibliography[1]{
  \OLDthebibliography{#1}
  \setlength{\itemsep}{0.375\baselineskip}
}

\hypersetup{
    linkbordercolor=bricks!75,
    citebordercolor=darts!75,
    linkcolor=black,
    urlbordercolor=fista!75,
    citecolor=black,
    pdftitle={Onboard DQG JGCD},
}
\title{Onboard Dual Quaternion Guidance for Rocket Landing}

\author{Abhinav G.\ Kamath$^*$, Taewan Kim$^*$, Skye Mceowen\footnote{Ph.D.\ Candidate, William E.\ Boeing Department of Aeronautics \& Astronautics; \texttt{\{agkamath,\,twankim,\,skye95\}@uw.edu}}, Mehran Mesbahi$^\dagger$, and \behcet\footnote{Professor, William E.\ Boeing Department of Aeronautics \& Astronautics; AIAA Fellow; \texttt{\{mesbahi,\,behcet\}@uw.edu}}}
\affil{University of Washington, Seattle, WA 98195, USA}

\author{Javier A. Doll\footnote{Senior Guidance Engineer; \texttt{javier.a.doll@nasa.gov}}}
\affil{Draper Laboratory, Houston, TX 77058, USA}

\author{Purnanand Elango\footnote{Research Scientist (Ph.D.\ Candidate at UW during the development of this work); \texttt{elango@merl.com}}}
\affil{Mitsubishi Electric Research Laboratories, Cambridge, MA 02139, USA}

\author{Yue Yu\footnote{Assistant Professor, Department of Aerospace Engineering and Mechanics; \texttt{yuey@umn.edu}}}
\affil{University of Minnesota Twin Cities, Minneapolis, MN 55455, USA}

\author{Taylor P. Reynolds\footnote{Senior Applied Scientist; \texttt{tayreyno@amazon.com}}}
\affil{Amazon Prime Air, Seattle, WA 98108, USA}

\author{Gavin F.\ Mendeck\footnote{Guidance, Navigation, and Control Flight Software Lead (SPLICE); \texttt{gavin.f.mendeck@nasa.gov}} and John M.\ Carson III\footnote{Technical Integration Manager – Precision Landing, NASA STMD; AIAA Fellow; \texttt{john.m.carson@nasa.gov}}}
\affil{NASA Johnson Space Center, Houston, TX 77058, USA}

\begin{document}

\maketitle

\begin{abstract}
\small
\vspace{0em}
The dual quaternion guidance (DQG) algorithm was selected as the candidate six-degree-of-freedom (6-DoF) powered-descent guidance algorithm for NASA's Safe and Precise Landing -- Integrated Capabilities Evolution (SPLICE) project. DQG is capable of handling state-triggered constraints that are of utmost importance in terms of enabling technologies such as terrain relative navigation (TRN). In this work, we develop a custom solver for DQG to enable onboard implementation for future rocket landing missions. We describe the design and implementation of a real-time-capable optimization framework, called sequential conic optimization (SeCO), that blends together sequential convex programming and first-order conic optimization to solve difficult nonconvex trajectory optimization problems, such as DQG, in real-time. This framework is entirely devoid of matrix factorizations/inversions, making it suitable for real-time applications. A key feature of SeCO is that it leverages a first-order primal-dual conic optimization solver, based on the proportional-integral projected gradient method (PIPG), that combines ideas pertaining to projected gradient descent and proportional-integral feedback of constraint violation. Unlike other conic optimization solvers, PIPG effectively exploits the sparsity structure and geometry of the constraints, avoids expensive equation-solving, and is suitable for both real-time and large-scale applications. We describe the implementation of this solver, develop customizable first-order methods, and leverage convergence-accelerating strategies such as warm-starting and extrapolation, to solve the nonconvex DQG optimal control problem in real-time. We show that the DQG-customized subproblem solver is able to solve the problem significantly faster than other state-of-the-art convex optimization solvers. Finally, in preparation for an upcoming closed-loop flight test campaign, we test our custom solver onboard the NASA SPLICE Descent and Landing Computer (DLC) in a hardware-in-the-loop setting. We observe that our algorithm is significantly faster than previously reported solve-times using the flight-tested interior point method (IPM)-based subproblem solver, BSOCP. Furthermore, our custom solver meets (and exceeds) NASA's autonomous precision rocket-landing guidance update-rate requirements for the first time, thus demonstrating the viability of SeCO for real-time, mission-critical applications onboard computationally-constrained flight hardware.
\vspace{0em}
\end{abstract}

\section{Introduction}

With robotic and human missions to the Moon and Mars on the horizon, there has been an increased interest in guidance, navigation, and control (GNC) technologies for precision landing \cite{smith2020artemis, chavers2019nasa, chavers2020nasa, petersen2020apollo, musk2017making, muirhead2020mars}. Precision landing and hazard avoidance (PL\&HA) have been deemed high-priority capabilities by NASA to facilitate missions of exploration to celestial bodies in the solar system \cite{carson2019splice}. Critical to achieving PL\&HA is powered-descent guidance (PDG), which refers to the generation of a feedforward control profile and the corresponding reference state trajectories for powered-descent and landing.

Historically, missions to the Moon and Mars—such as the Apollo program, and the Mars Science Laboratory (MSL) and Mars 2020 missions, respectively—made use of polynomial guidance for the powered-descent and landing phase \citep{klumpp1974apollo, san2013development, casoliva2021reconstructed}. While these missions were highly successful, the resulting descent and landing trajectories were not propellant optimal. Further, the landing dispersion ellipses for the Apollo missions, spanning kilometers \cite{quaide1969geology}, were too large for precision landing, which requires safely landing spacecraft within 100 meters of the target landing site \cite{carson2019splice}.

The advent of convex optimization methods at the turn of the century, along with pivotal results on the lossless convexification ({\lcvx}) of certain nonconvex constraints, led to the development of a real-time-capable PDG algorithm with strong guarantees \cite{acikmese2007convex, accikmecse2011lossless, accikmecse2013lossless}. This algorithm was successfully flight-tested onboard a terrestrial rocket-powered landing testbed, demonstrating the capability of generating propellant-optimal large divert trajectories in real-time \cite{acikmese2013flight, scharf2017implementation}. It has also been implemented in the context of upcoming lunar landing missions \cite{berning2023lossless, shaffer2024implementation}. The problem formulation incorporates a point-mass vehicle model and three-degree-of-freedom (3-DoF) translation dynamics. While this algorithm relies on a direct optimal control method with either a zero-order hold (ZOH) \cite{behcet2007jgcd} or a first-order hold (FOH) \cite{accikmecse2013lossless} control parameterization, it has also been extended to work within a direct pseudospectral optimal control framework \cite{sagliano2018pseudospectral, sagliano2019generalized}, an indirect optimal control framework \cite{lu2018propellant}, and a hybrid direct-indirect optimal control framework \cite{spada2023direct}.

While it has been demonstrated that a 3-DoF guidance algorithm is sufficient to successfully guide and land an inherently 6-DoF vehicle (by using the thrust vector profile as a surrogate for the reference attitude of the vehicle) \cite{acikmese2013flight, scharf2017implementation}, the lack of explicit attitude states in the 3-DoF model makes it challenging to impose constraints on the attitude of the vehicle. This challenge has motivated the development of planar landing guidance algorithms with attitude modeling \cite{reynolds2020optimal, reynolds2020real} and explicit 6-DoF powered-descent guidance algorithms \cite{szmuk2017successive, szmuk2018successive, spada2023successive, chari2024fast, sagliano2024six, elango2025continuous}, both of which are inherently nonconvex, necessitating the use of nonconvex optimization methods such as successive convexification, which falls under the umbrella of sequential convex programming (SCP) \cite{malyuta2021convex, elango2025successive}.

A key challenge in landing guidance algorithms is constraining the flight envelope in a manner that couples the translation and attitude states. The resulting state constraints are critical for technologies such as terrain relative navigation (TRN) and hazard detection and avoidance (HDA). Moreover, they are often required to be imposed conditionally. Recent work on 6-DoF guidance addresses this problem by a combination of coupled-constraint modeling \cite{lee2015optimal, lee2017constrained, hayner2025continuous} and the formulation of \emph{state-triggered constraints}, (which are then encoded in a continuous optimization framework) \cite{szmuk2019successive, reynolds2019state, szmuk2020successive, reynolds2020dual, kamath2023real, buckner2024constrained, uzun2025sequential, kim2025six}, and the use of either piecewise-affine model predictive control (PWA-MPC) or SCP for trajectory optimization. These are powerful algorithms capable of handling challenging constraints, and the SCP-based methods, specifically, have been shown to be amenable to real-time implementation \cite{szmuk2020successive, reynolds2020dual}.

One such state-of-the-art algorithm is the dual quaternion guidance (\dqg) algorithm for rocket landing, which was originally presented in \citep{lee2017constrained, reynolds2019state, reynolds2020dual}; the algorithm in \citep{reynolds2020dual} was chosen as the candidate 6-DoF powered-descent guidance algorithm for NASA's Safe and Precise Landing – Integrated Capabilities Evolution (SPLICE) project, and has been open-loop flight-tested on the Blue Origin New Shepard suborbital rocket onboard the Descent and Landing Computer (DLC) \cite{rutishauser2021nasa, fritz2022post}.

The implementation of the algorithm for these flight tests made use of a customized version of an interior point method (IPM)-based convex subproblem solver called {\bsocp}, used in conjunction with a parser interface called {\cprs} \cite{dueri2014automated, dueri2017customized}. For the terrestrial test flights, {\dqg} was solved in under 3 seconds (in an open-loop, i.e., the generated solution was not utilized by the rocket), although it was executed at 0.2 Hz (once every 5 seconds) onboard the DLC \cite{fritz2022post}. Customization of {\bsocp} for this application led to a very large footprint source code—around 600,000 lines of C code. As noted in \cite{strohl2022implementation}, that formulation of {\dqg} was later adapted to lunar powered-descent, but was implemented with the generic (uncustomized) version of {\bsocp} instead—this implementation took over 11 seconds to execute with an optimization horizon length of 20, and close to 6 seconds to execute with an optimization horizon length of 10, making it prohibitively slow in terms of meeting NASA's guidance update-rate requirements for PL\&HA in its current state \cite{doll2025hardware}.

In order to tackle the challenges of execution speed and code footprint, first-order optimization algorithms can be used instead of (second-order) IPMs. First-order algorithms typically rely on simple linear algebra operations like matrix-vector multiplications and computation of vector norms at each iteration. Unlike second-order methods, they can avoid factorization of larger matrices and can be warm-started easily. In addition to these benefits, the recently introduced first-order algorithm, the proportional-integral projected gradient ({\pipg}) method \cite{yu2020proportional, yu2021proportional}, is capable of exploiting the structure of trajectory optimization problems to completely avoid operations on large sparse matrices. As a result, {\pipg} is readily suitable for onboard, resource-constrained applications. 

Recent work on the design and implementation of a customized {\pipg} solver for the 3-DoF {\lcvx} algorithm demonstrated the capabilities of {\pipg} in terms of computational efficiency, code footprint, and its viability for easy verification and validation \cite{elango2022customized}. Moreover, features such as warm-starting and extrapolation \cite{yu2022extrapolated} for boosting the practical convergence rate have further enhanced capabilities of {\pipg} as a conic optimization solver for use within sequential convex programming (SCP) algorithms \cite{SCPToolboxCSM2022}. In fact, SCP algorithms can be specialized to harness the approach taken in {\pipg} to solve conic optimization problems, as shown in \cite{kamath2023real}, with the sequential conic optimization (\seco) framework, which is entirely devoid of matrix factorizations and inversions.

In this work, we demonstrate a real-time-capable implementation of {\seco} for solving the nonconvex optimal control problem in {\dqg}. Sections \ref{sec:ocp}, \ref{sec:dynamics}, \ref{sec:seco}, and \ref{sec:solver} describe the problem formulation and optimization algorithms, Section \ref{sec:custom} delves into solver customization, and Section \ref{sec:results} provides the numerical results.

This work is a direct extension of our initial conference paper \cite{kamath2023customized}, addressing both the tasks listed under the planned future work. Specifically, we: (i) present an extended treatment of the preconditioning procedure, borrowing from \cite{kamath2025optimal} (the work in which was motivated by the open problem in the initial conference paper) to automatically tune the solver parameter—we provide a customized version of the algorithm for the same in this extension; and, (ii) we provide hardware-in-the-loop test results from onboard the SPLICE Descent and Landing Computer (DLC), with realistic mission parameters, for an upcoming closed-loop rocket landing flight test campaign \citep{Mendeck_SPLICE_2023, Mendeck_SPLICE_2024}, that were originally presented in another conference paper \cite{doll2025hardware}. We note that there have been recent advancements in the field, such as the continuous-time successive convexification (\ctscvx{}) framework that ensures continuous-time constraint satisfaction \cite{elango2025successive, uzun2024successive}, and the auto-tuned primal-dual successive convexification (\autoscvx{}) that enables the automatic tuning of SCP hyperparameters \cite{mceowen2025auto}; these new techniques will not be discussed in this work.

Many of the following fundamental mathematical definitions and formalisms can be found in \cite{reynolds2020dual, reynolds2020computational, kamath2023real}, and are presented here for the sake of completeness. The operations pertaining to quaternion and dual quaternion algebra can be found in \cite{lee2012dual, lee2015optimal, lee2017constrained, reynolds2019state, reynolds2020dual, reynolds2020computational}, and are documented in the \hyperref[appendix]{Appendix}.
\section{Optimal Control Problem Formulation}\label{sec:ocp}

In this section, we present the continuous-time optimal control formalism for
the 6-DoF rocket landing guidance problem, including the equations of
motion and pertinent constraints on the state and
control.

\subsection{Equations of Motion}

One of the defining characteristics of {\dqg} is the representation of the 6-DoF equations of motion using unit dual quaternions, which yield an elegant parameterization of the dynamics by naturally coupling the translational and rotational states and enabling the representation of certain key operational constraints, such as the line-of-sight constraint, as convex constraints (in theory) \cite{lee2015optimal, reynolds2020dual}. However, the use of the dual quaternion parameterization is ultimately a design choice, and other parameterizations can be adopted as well, such as Cartesian coordinates for the translational states and unit quaternions for attitude \cite{szmuk2020successive}. The interested reader is referred to \citep{lee2012dual, lee2017constrained,  reynolds2018coupled, reynolds2020dual} for detailed descriptions of parameterizing rigid body dynamics via unit dual quaternions.

\subsubsection{States}

The state vector is 15-dimensional, and consists of mass, $m$, the 8-dimensional unit dual quaternion that couples translation and attitude, $\dq{q}$, and the (reduced-order) 6-dimensional dual velocity, $\dq{\omega}$, as shown in Equations \eqref{eq:state_vector}. The vector $\widetilde{\dq{\omega}}$ represents the 8-dimensional dual velocity, in which the fourth and eighth terms are zero; $\q{q}$ is the attitude (unit) quaternion. In this work, we adopt the scalar-last convention to represent quaternions. See the \hyperref[appendix]{Appendix} for definitions of unit quaternions and unit dual quaternions. The subscripts $\I$ and $\B$ denote that the quantity in question is expressed in the inertial frame or the body frame, respectively. Further, $r \inthree$ is the position, $\omega \inthree$ is the angular velocity, and $v \inthree$ is the velocity; the aforementioned representations are summarized as follows:
\begin{subequations}
\allowdisplaybreaks
\begin{align}
    m &\inone\\
    \dq{q}\ &{\defeq} \begin{pmatrix}\q{q}\\\dfrac{1}{2}\!\begin{pmatrix}\q{r}_{\I}\\ 0\end{pmatrix}\otimes\q{q}\end{pmatrix} = \begin{pmatrix}\q{q}\\\dfrac{1}{2}\q{q}\otimes\begin{pmatrix}\q{r}_{\B}\\ 0\end{pmatrix}\end{pmatrix} \in\mathbb{R}^{8}\\
    \widetilde{\dq{\omega}}\ &{\defeq} \begin{pmatrix}\begin{pmatrix}\q{\omega}_{\B}\\ 0\end{pmatrix}\\[.2in]q^{*} \otimes \begin{pmatrix}\q{v}_{\I}\\ 0\end{pmatrix} \otimes q\end{pmatrix} = \begin{pmatrix}\begin{pmatrix}\q{\omega}_{\B}\\ 0\end{pmatrix}\\[.2in]\begin{pmatrix}\q{v}_{\B}\\ 0\end{pmatrix}\end{pmatrix} \ineight\\
    \dq{\omega}\ &{\defeq} \begin{pmatrix}\q{\omega}_{\B}\\\q{v}_{\B}\end{pmatrix} \insix\\
    x\ &{\defeq} \begin{pmatrix}m\\\dq{q}\\\dq{\omega}\end{pmatrix}\infifteen
    \label{eq:state_vector}
\end{align}
\end{subequations}

\subsubsection{Controls}

The control input vector is 6-dimensional, with three parameters describing the thrust vector: the thrust magnitude, $T \inone$, the gimbal deflection angle, $\delta \inone$, and the gimbal azimuth angle, $\phi \inone$, and a 3-dimensional body torque vector, $\tau \inthree$, as shown in Equations \eqref{eq:control_input_vector}. The thrust vector is effected by means of a gimbaled main engine and the torque input is assumed to be effected by means of reaction control system (RCS) thrusters.
\begin{figure}[H]
    \begin{mybox}
        \centering
        \includegraphics{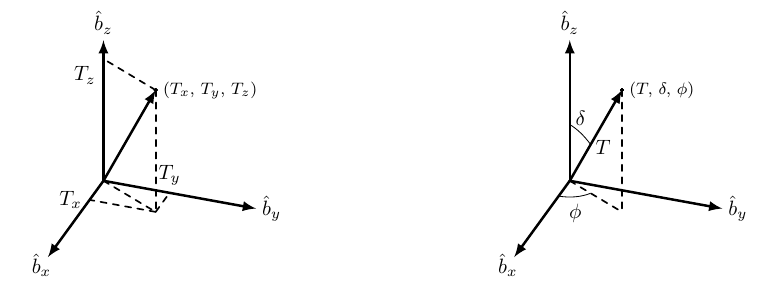}
    \end{mybox}
    \caption{Parameterization of the thrust vector (expressed in the body frame) in terms of Cartesian coordinates (left) and spherical coordinates (right). In this work, we use the latter.}
    \label{fig:control_parameterizations}
\end{figure}
\vspace{-1em}
We choose to parameterize the thrust vector in terms of spherical coordinates, as shown in Figure \ref{fig:control_parameterizations}, for the following reasons:
(1) all the control constraints in {\dqg} become naturally convex, and hence, in combination with the first-order hold (FOH) parameterization (described in Subsection \ref{subsec:linearization}), intersample satisfaction of the control constraints is guaranteed; (2) the control rate constraints can be imposed exactly; and, (3) with a mild assumption, both the magnitude and rate constraints (throttle and gimbaling) can be combined and made projection-friendly, which is beneficial in terms of implementing the solver, as described in Subsection \ref{subsec:constraint_classification}. Further, with the spherical coordinate parameterization, we note that the thrust magnitude solution possesses a piecewise-affine profile, which will not be the case if the Cartesian coordinate parameterization is adopted; these parameterizations are summarized as follows:
\begin{subequations}
\allowdisplaybreaks
    \begin{align}
        \mathcal{T}_{\B} = \begin{pmatrix} T_{x} \\ T_{y} \\ T_{z} \\ \tau_{x} \\ \tau_{y} \\ \tau_{z} \end{pmatrix} &= \begin{pmatrix}T \sin \delta \cos \phi \\ T \sin \delta \sin \phi \\ T \cos \delta \\ \tau_{x} \\ \tau_{y} \\ \tau_{z} \end{pmatrix} \insix \label{eq:trig_nonlinear}\\
        u &= \begin{pmatrix}T \\ \delta \\ \phi \\ \tau_{x} \\ \tau_{y} \\ \tau_{z} \end{pmatrix} \insix
    \end{align}
    \label{eq:control_input_vector}
\end{subequations}
\begin{center}
  \begin{tabular}{l}
        $T$: thrust magnitude\\
        $\delta$: gimbal deflection angle defined from the body vertical ($z$) axis\\
        $\phi$: gimbal azimuth angle defined from the body $x$-axis\\
        $\tau_{x}$, $\tau_{y}$, $\tau_{z}$: body torque inputs
  \end{tabular}
 \end{center}
Here, $\mathcal{T}_{\B}$ is the wrench vector expressed in the body frame \citep{szmuk2019successive}, and $u$ is the control input vector. The convexity and simplicity of the resulting control constraints come at the cost of additional trigonometric nonlinearities in the dynamics, as shown in Equation \eqref{eq:trig_nonlinear}.

\subsubsection{Mass-Depletion}
In addition to accounting for thrust due to the main engine, we consider the effect of thrust due to the RCS thrusters on mass-depletion. In order to do so, we assume that two diagonally opposite RCS thrusters fire at any given instant to achieve the desired net torque, such that the thrust due to each thruster is orthogonal to the body vertical ($z$) axis.

Assuming the mass-center of the vehicle is equidistant from the top-mounted RCS thrusters and the bottom-mounted RCS thrusters/gimbaled main engine, the mass-depletion dynamics can be given by Equation \eqref{eq:mass_depletion}:
\begin{align}
    \dot{m}(t) = -\!\left(\alpha_{\textsc{me}}\,T(t) + \alpha_{\textsc{rcs}}\frac{\norm{\tau(t)}_{2}}{l_{\textsc{cm}}}\right)
    \label{eq:mass_depletion}
\end{align}
Here, $\alpha_{\textsc{me}} \in \mathbb{R}_{+}$ and $\alpha_{\textsc{rcs}} \in \mathbb{R}_{+}$ are the thrust-specific fuel consumption (TSFC) parameters for the main engine and an RCS thruster, respectively, and $l_{\textsc{cm}} \in \mathbb{R}_{++}$ is the length of the moment-arm of the vehicle, i.e., the distance between the mass-center of the vehicle and the bottom-mounted RCS thrusters/gimbaled main engine.

\subsubsection{Kinematics}
This dual quaternion kinematic equation requires the use of the 8-dimensional dual velocity vector, where the fourth and eighth terms are zero, and is given by Equation \eqref{eq:kinematics}:
\begin{align}
    \dot{\dq{q}}(t) = \frac{1}{2}\dq{q}(t)\otimes\widetilde{\dq{\omega}}(t)
    \label{eq:kinematics}
\end{align}

\subsubsection{Dynamics}
For the dynamics, we assume that the moment of inertia is a linear function of the vehicle mass (as opposed to assuming a constant value as in \citep{reynolds2020dual}), and thereby account for the effect of mass-depletion on the attitude of the vehicle. The moment of inertia can be assumed to be an affine function of the vehicle mass as well, if required \cite{reynolds2019state}. The dynamics can be given by Equation \eqref{eq:dynamics}:
\begin{align}
    \dot{\dq{\omega}}(t) = \dq{J}(t)^{-1}\!\left[\begin{pmatrix}0_{3\times3} & -\q{\omega}_{\B}^{\times}\!\!\>(t)\\[1ex]-\q{\omega}_{\B}^{\times}\!\!\>(t) & 0_{3\times3}\end{pmatrix}\dq{J}(t)\!\>\dq{\omega}(t) + \q{\Phi}\!\>\mathcal{T}_{\B}(t) + m(t)\!\>\dq{g}_{\B}(t)\right]
    \label{eq:dynamics}
\end{align}
where
\[
  \dq{J}(t) \defeq
  \left(\begin{array}{c|c}
  0_{3\times3} &
  m(t)\!\>I_{3}\\[1ex]
  \hline\\[-2ex]
  m(t)\!\>\q{J} &
  0_{3\times3}
  \end{array}\right)_{6\times6}
  \qquad
  \q{\Phi} \defeq
  \left(\begin{array}{c|c}
  ~I_{3}~&
  0_{3\times3}\\[1ex]
  \hline\\[-2ex]
  l^{\times} &
  I_{3}
  \end{array}\right)_{6\times6}
  \qquad
  \dq{g}_{\B}(t) \defeq
  \left(\begin{array}{c}
  \q{g}_{\B}(t)\\[1ex]
  \hline\\[-2ex]
  0_{3\times1}
  \end{array}\right)_{6\times1}
\]
\begin{align*}
    \q{g}_{\I}\ &{\defeq}~\left[0,\,0,\,-g\right]^{\top}\\[5ex]
    \begin{pmatrix}\q{g}_{\B}(t)\\0\end{pmatrix} &= \q{q}^{*}\!\!\>(t)\otimes\begin{pmatrix}\q{g}_{\I}\\0\end{pmatrix}\otimes\q{q}(t)
\end{align*}
where $J \in \mathbb{S}^{3}_{++}$ is the inertia tensor of the vehicle about its mass-center, $l = \left[0, 0, -l_{\textsc{cm}}\right]^{\top}$ is the body-fixed moment-arm vector, and $g \in \mathbb{R}_{+}$ is the acceleration due to gravity at the celestial body under consideration.
\subsection{Control Constraints}
All components of the thrust vector are bounded, as shown in Equations \eqref{eq:control_magnitude}, where $T_{\min} \in \mathbb{R}_{++}$ and $T_{\max} \in \mathbb{R}_{++}$ are the lower and upper bounds on the thrust magnitude, respectively, and $\delta_{\max} \in \mathbb{R}_{++}$ is the upper bound on the gimbal deflection angle.  Further, they are rate-limited, as shown in Equations \eqref{eq:control_rate}, where $\dot{T}_{\max} \in \mathbb{R}_{++}$, $\dot{\delta}_{\max} \in \mathbb{R}_{++}$, and $\dot{\phi}_{\max} \in \mathbb{R}_{++}$ are the rate-limits on the thrust magnitude, the gimbal deflection angle, and the gimbal azimuth angle, respectively. An upper bound is levied on the magnitude of the body torque input about each body axis, as shown in \eqref{eq:torque_bounds}. We emphasize that every single control constraint is naturally convex, owing to the spherical coordinate parameterization.
\subsubsection{Thrust Vector Bounds}
\vspace{-1em}
\begin{subequations}
\begin{align}
    T_{\min} &\le T(t) \le T_{\max} \\
    0 &\le \delta(t) \le \delta_{\max} \\
    0 &\le \phi(t) \le 2\pi
\end{align}
\label{eq:control_magnitude}
\end{subequations}
\vspace{-1em}
\subsubsection{Thrust Vector Rate-Limits}
\vspace{-1em}
\begin{subequations}
\begin{align}
     \abs{\dot{T}(t)} &\le \dot{T}_{\max} \\
     \abs{\dot{\delta}(t)} &\le \dot{\delta}_{\max} \\
     \abs{\dot{\phi}(t)} &\le \dot{\phi}_{\max}
\end{align}
\label{eq:control_rate}
\end{subequations}
\vspace{-1em}
\subsubsection{Torque Bounds}
\vspace{-1em}
\begin{align}
    \norm{\tau(t)}_{\infty} \le \tau_{\max}
    \label{eq:torque_bounds}
\end{align}
\subsection{State Constraints}
We classify state constraints into global state constraints, state-triggered constraints, initial condition constraints, and terminal condition constraints, which we describe in this subsection.
\subsubsection{Global State Constraints}
Global state constraints involve constraints that are imposed on the state over the entire time-horizon. These constraints include a maximum tilt constraint, a maximum angular body rate constraint, a maximum speed constraint, and a minimum altitude constraint, as shown in Equations \eqref{eq:global}, respectively:
\begin{subequations}
\begin{align}
    \forall t \in \left[0, t_{f}\right), \qquad\qquad& \nonumber\\
    \norm{\dq{q}^{[1:2]}(t)}_{2} &\le \sin \frac{\theta_{\max}}{2} \\
    \norm{\dq{\omega}^{[1:3]}(t)}_{\infty} &\le \omega_{\max} \\
    \norm{\dq{\omega}^{[4:6]}(t)}_{2} &\le v_{\max} \\
    \dq{q}(t)^{\top} M_{g}\,\dq{q}(t) &\ge h_{\min} \label{eq:min_alt_ncvx}
\end{align}
\label{eq:global}
\end{subequations}
where
\begin{align*}
    M_g \defeq \begin{pmatrix}
    \mathbf{0}_{4 \times 4} & {\left[\widetilde{z}_{\mathcal{I}}\right]_{\otimes}^{\top}} \\
    {\left[\widetilde{z}_{\mathcal{I}}\right]_{\otimes}} & \mathbf{0}_{4 \times 4}
    \end{pmatrix}
\end{align*}
Here, $t_{f} \in \mathbb{R}_{+}$ is the time-of-flight, $\theta_{\max} \in \mathbb{R}_{+}$ is the maximum tilt angle from the inertial vertical ($z$) axis, $\omega_{\max} \in \mathbb{R}_{+}$ is the maximum angular speed about any body axis, $v_{\max} \in \mathbb{R}_{+}$ is the maximum speed, $h_{\min} \in \mathbb{R}_{+}$ is the minimum altitude, $z_{\I} \defeq [0, 0, 1]^{\top}$, and $\widetilde{z}_{\I} \defeq [z_{\I}^{\top}, 0]^{\top}$ (pure quaternion).
\subsubsection{State-Triggered Constraints}
State-triggered constraints (STCs) include the constraints that are to be activated only when the vehicle is within the prescribed trigger window. Here, we consider slant-range-based triggering, and tightly constrain the body tilt angle, angular body rates, maximum speed, and the maximum line-of-sight angle to the target landing site---which is assumed to be at the origin, without loss of generality---as shown in Equations \eqref{eq:state_triggered}, respectively:
\begin{subequations}
\allowdisplaybreaks
\begin{align}
\forall t \ni \rho_{\min} \le \norm{2\,\dq{q}^{[5:8]}(t)}_{2} \le \rho_{\max}, \qquad& \nonumber\\
    \norm{\dq{q}^{[1:2]}(t)}_{2} &\le \sin \frac{\theta_{\textsc{stc}_{\max}}}{2} \\
    \norm{\dq{\omega}^{[1:3]}(t)}_{\infty} &\le \omega_{\textsc{stc}_{\max}} \\
    \norm{\dq{\omega}^{[4:6]}(t)}_{2} &\le v_{\textsc{stc}_{\max}} \\
    \dq{q}(t)^{\top} M_{l}\,\dq{q}(t) &+ \norm{2\,\dq{q}^{[5:8]}(t)}_{2} \cos \mu_{\textsc{stc}_{\max}} \le 0 \label{eq:LoS_ncvx}
\end{align}
\label{eq:state_triggered}%
\end{subequations}
where
\begin{align*}
    M_l \defeq \begin{pmatrix}
    \mathbf{0}_{4 \times 4} & {\left[\widetilde{p}_{\mathcal{B}}\right]_{\otimes}^{* \top}} \\
    {\left[\widetilde{p}_{\mathcal{B}}\right]_{\otimes}^*} & \mathbf{0}_{4 \times 4}
    \end{pmatrix}
\end{align*}
Here, $\rho_{\max} \in \mathbb{R}_{++}$ and $\rho_{\min} \in \mathbb{R}_{++}$ are the maximum (activation) and minimum (deactivation) trigger distances from the target landing site, respectively; $\theta_{\textsc{stc}_{\max}}$, $\omega_{\textsc{stc}_{\max}}$, and $v_{\textsc{stc}_{\max}}$ are assumed to be smaller than their counterparts in Equations \eqref{eq:global}; $\mu_{\textsc{stc}_{\max}} \in \mathbb{R}_{+}$ is the maximum line-of-sight angle; and $\widetilde{p}_{\mathcal{B}} \defeq [p_{\mathcal{B}}^{\top}, 0]^{\top}$ (pure quaternion), where $p_{\B} \inthree$ is a unit vector in the body frame that represents the body-fixed sensor-pointing direction. We choose not to impose the minimum altitude constraint in the trigger window, with the observation that the simultaneous satisfaction of the maximum tilt and maximum line-of-sight angle constraints implicitly precludes subsurface solutions if $\theta_{\textsc{stc}_{\max}} \le \frac{\pi}{2} - \mu_{\textsc{stc}_{\max}} - \gamma_{\text{boresight}}$, where $\gamma_{\text{boresight}}$ is the angle made by the body-fixed sensor with the body vertical ($z$) axis, i.e., the (acute) angle between $p_{\B}$ and $\hat{b}_{z}$.

These constraints are imposed to enable accurate scans of the potential landing site during descent using the hazard detection LiDAR (HDL) \cite{restrepo2019nasa, restrepo2020next}, for instance, and to initiate diverts if necessary.

\subsubsection{Initial Conditions}
The initial condition constraints are given by Equations \eqref{eq:initial_conditions}:
\begin{subequations}
\begin{align}
    m(0) &= m_{i} \\
    \dq{q}(0) &= \begin{pmatrix}\q{q}_{i} \\ \dfrac{1}{2}\!\begin{pmatrix}\q{r}_{\I_{i}}\\ 0\end{pmatrix}\otimes\q{q}_{i}\end{pmatrix} \\
    \dq{\omega}(0) &= \begin{pmatrix}\q{\omega}_{\B_{i}}\\q_{i}^{*} \otimes \begin{pmatrix}\q{v}_{\I_{i}}\\ 0\end{pmatrix} \otimes q_{i}\end{pmatrix}
\end{align}
\label{eq:initial_conditions}%
\end{subequations}
where $m_{i} \in \mathbb{R}_{++}$ is the initial mass of the vehicle, $q_{i} \infour_{u}$ is the initial attitude quaternion, $r_{\I_{i}} \inthree$ is the initial position expressed in the inertial frame, $\omega_{\B_{i}} \inthree$ is the initial angular velocity expressed in the body frame, and $v_{\I_{i}} \inthree$ is the initial velocity expressed in the inertial frame.
\subsubsection{Terminal Conditions}
The terminal conditions subsume the following details—at the final time $t_{f}$: (1) the vehicle is upright (zero pitch and yaw); (2) the roll is free; (3) the angular body rates are zero; and, (4) the (inertial) horizontal components of velocity are zero. By infusing these details into the expressions rather than treating them as problem parameters, the terminal condition constraints are rendered convex. Here, $m_{f} \in \mathbb{R}_{++}$ is the final mass of the vehicle, $r_{\I_{f}} \inthree$ is the final position expressed in the inertial frame, and $v_{z_{\I_{f}}}\inone$ is the final velocity along the inertial vertical ($z$) axis. If desired, the entire final dual quaternion can be fixed as well. The terminal condition constraints are given by Equations \eqref{eq:terminal_conditions}:
\begin{subequations}
\allowdisplaybreaks
\begin{align}
    m(t_{f}) &\ge m_{f} \label{eq:mass_lower_bound}\\
    \dq{q}^{[1:2]}(t_{f}) &= \dq{0}_{2 \times 1} \\
    \begin{bmatrix} \left(\frac{1}{2}[r_{\mathcal{I}_f}]_{\otimes}\!\begin{pmatrix}\mathbf{0}_{2\times2}\\ I_{2}\end{pmatrix}\right) & -I_{4} \end{bmatrix}_{\!4 \times 6}\!\dq{q}^{[3:8]}(t_{f}) &= \dq{0}_{4 \times 1}\\
    \widetilde{\dq{\omega}}(t_{f}) &= \begin{pmatrix}\dq{0}_{4 \times 1}\\
    \begin{pmatrix} \dq{0}_{2 \times 1} \\ q^{[3:4]}(t_{f}) \end{pmatrix}^{\!*} \otimes \begin{pmatrix} \dq{0}_{2 \times 1} \\ v_{z_{\I_{f}}} \\ 0 \end{pmatrix} \otimes \begin{pmatrix} \dq{0}_{2 \times 1} \\ q^{[3:4]}(t_{f}) \end{pmatrix}
    \end{pmatrix} = \begin{pmatrix}
    \dq{0}_{4 \times 1} \\ \dq{0}_{2 \times 1} \\ v_{z_{\I_{f}}} \\ 0
    \end{pmatrix} \\
    \therefore \dq{\omega}(t_{f}) &= \begin{pmatrix}
    \dq{0}_{5 \times 1} \\ v_{z_{\I_{f}}}
    \end{pmatrix}
\end{align}
\label{eq:terminal_conditions}%
\end{subequations}
\begin{mybox}
    \vspace{-0.75em}
    \subsection{The Continuous-Time Nonconvex Optimal Control Problem}\label{subsec:cont_time_ocp}
        \begin{equation*}
        \small
            \begin{aligned}
            \underset{t_{f},\, u(t)}{\text{minimize}} &
            && -m(t_{f}) \\
            \text{subject to} &
            && \forall t \in [0, t_{f}) \\
            \fbox{\text{Dynamics}} &&& {\f{\dot{x}}{t}} = \f{f}{{t, {\f{x}{t}}}, {\f{u}{t}}} \\
            \fbox{\text{Control constraints}} & && T_{\min} \le T(t) \le T_{\max} \\
            &&& 0 \le \delta(t) \le \delta_{\max} \\
            &&& 0 \le \phi(t) \le 2\pi \\
            &&& \abs{\dot{T}(t)} \le \dot{T}_{\max} \\
            &&& \abs{\dot{\delta}(t)} \le \dot{\delta}_{\max} \\
            &&& \abs{\dot{\phi}(t)} \le \dot{\phi}_{\max} \\
            &&& \norm{\tau(t)}_{\infty} \le \tau_{\max} \\
            \fbox{\text{Global state constraints}} &&& 
            \norm{\dq{q}^{[1:2]}(t)}_{2} \le \sin \frac{\theta_{\max}}{2} \hphantom{\fbox{\text{Global state constraints~~~~~~~~~~~~}}} \\
            &&& \norm{\dq{\omega}^{[1:3]}(t)}_{\infty} \le \omega_{\max} \\
            &&& \norm{\dq{\omega}^{[4:6]}(t)}_{2} \le v_{\max} \\
            &&& \dq{q}(t)^{\top} M_{g}\,\dq{q}(t) \ge h_{\min} \\
            \fbox{\text{State-triggered constraints}} &&& \norm{\dq{q}^{[1:2]}(t)}_{2} \le \sin \frac{\theta_{\textsc{stc}_{\max}}}{2} \\
            \forall t \ni \rho_{\min} \le \norm{2\,\dq{q}^{[5:8]}(t)}_{2} \le \rho_{\max} &&& \norm{\dq{\omega}^{[1:3]}(t)}_{\infty} \le \omega_{\textsc{stc}_{\max}} \\
            &&& \norm{\dq{\omega}^{[4:6]}(t)}_{2} \le v_{\textsc{stc}_{\max}} \\
            &&& \dq{q}(t)^{\top} M_{l}\,\dq{q}(t) + \norm{2\,\dq{q}^{[5:8]}(t)}_{2} \cos \mu_{\textsc{stc}_{\max}} \le 0 \\
            \fbox{\text{Initial conditions}} &&& m(0) = m_{i} \\
            &&& \dq{q}(0) = \begin{pmatrix}\q{q}_{i} \\ \dfrac{1}{2}\!\begin{pmatrix}\q{r}_{\I_{i}}\\ 0\end{pmatrix}\otimes\q{q}_{i}\end{pmatrix} \\
            &&& \dq{\omega}(0) = \begin{pmatrix}\q{\omega}_{\B_{i}}\\q_{i}^{*} \otimes \begin{pmatrix}\q{v}_{\I_{i}}\\ 0\end{pmatrix} \otimes q_{i}\end{pmatrix} \\
            \fbox{\text{Terminal conditions}} &&& m(t_{f}) \ge m_{f} \\
            &&& \dq{q}^{[1:2]}(t_{f}) = \dq{0}_{2 \times 1} \\
            &&& \begin{bmatrix} \left(\frac{1}{2}[r_{\mathcal{I}_f}]_{\otimes}\!\begin{pmatrix}\mathbf{0}_{2\times2} \\ I_{2}\end{pmatrix}\right) & -I_{4} \end{bmatrix}_{\!4 \times 6}\!\dq{q}^{[3:8]}(t_{f}) = \dq{0}_{4 \times 1} \\
            &&& \dq{\omega}(t_{f}) = \begin{pmatrix}
            \dq{0}_{5 \times 1} \\ v_{z_{\I_{f}}}
            \end{pmatrix}
            \end{aligned}
        \end{equation*}
        \vspace{0em}
\end{mybox}
\vspace{\baselineskip}
The continuous-time nonconvex optimal control problem is given by Subsection \ref{subsec:cont_time_ocp}, where we minimize propellant consumption (by maximizing the final mass of the vehicle).
\section{Transformation and Discretization of Dynamics}\label{sec:dynamics}

Our approach to treating the dynamics\footnote{Henceforth, we overload the term ``dynamics'' to encompass all the equations of motion.} closely follows the methods provided in \citep{szmuk2020successive, reynolds2020dual, reynolds2020real}, with some key distinctions, such as the inverse-free propagation step, as described later in this section. All of these approaches yield an \textit{exact} discretization of the linear time-varying (LTV) system under consideration---which means that the discrete-time trajectory exactly coincides with the continuous-time trajectory at the temporal nodes---and are analytically equivalent. In practice, however, the approach we propose herein leads to a much simpler implementation, without the need for any matrix factorizations/inversions, thus also making the implementation more numerically stable and reliable. Further, a very similar approach, albeit in a multi-phase guidance setting, is provided in \citep{kamath2023real}.

\subsection{Time-Dilation}

The original nonlinear dynamics over the entire time-horizon are given by Equation \eqref{eq:nonlinear}:
\begin{align}
    \f{\dot{x}}{t} = {\f{f}{t, {\f{x}{t}}, {\f{u}{t}}}}, \enskip t \in [0, t_{f}) \label{eq:nonlinear}
\end{align}
where, without loss of generality, we assume that the initial time is zero.
Now, we define an invertible linear map, $\tau:[0,t_f) \to [0,1)$, as shown in Equation \eqref{eq:dilation}:
\begin{align}
    {\f{\tau}{t}} = \dfrac{t}{t^{-}_{f}}, \enskip t \in [0, t_{f}) \label{eq:dilation}
\end{align}
where the negative-sign superscript denotes the left limit. This mapping is referred to as \textit{time-dilation}, as it dilates (normalizes) the wall-clock time-horizon to a chosen fixed interval—$[0, 1)$, in our case. Next, we invoke the chain-rule to obtain Equation \eqref{eq:chain}:
\begin{flalign}
{\f{\circdot{x}}{t}} &= \frac{d}{d\tau}\,
{\f{x}{t}} = \frac{d t}{d\tau}~\frac{d}{d t}\,
{\f{x}{t}} = \frac{d t}{d\tau}~{\f{\dot{x}}{t}} = \underbracket[0.14ex]{t^{-}_{f}}_{s}\,{\f{\dot{x}}{t}} \nonumber\\
&= s\,{\f{f}{t, {\f{x}{t}}, {\f{u}{t}}}} \defeq {\f{F}{t, {\f{x}{t}}, {\f{u}{t}}, s}} \label{eq:chain}
\end{flalign}
where $\circdot{\square}$ denotes the derivative operator with respect to the dilated time, $\tau$.
The multiplier in Equation \eqref{eq:chain}, $s \defeq t^{-}_{f} \in \R_{+}$, termed the \textit{dilation factor}, evaluates to the time-of-flight. Given Equations \eqref{eq:dilation} and \eqref{eq:chain}, $t$ can now be expressed as a function of $\tau$, i.e., $t(\tau) = s\,\tau$, $\tau \in [0, 1)$. Henceforth, we treat $\tau$ as the independent variable instead of $t$, and replace the temporal argument, $t = t(\tau)$, by $\tau$, for notational simplicity. Note that time-dilation transforms the free-final-time optimal control problem to an equivalent fixed-final-time optimal control problem.


\subsection{Linearization}\label{subsec:linearization}

We begin by considering nonlinear systems that can be represented as ordinary differential equations (ODEs), as given by Equation \eqref{eq:ode}:
\begin{align}
    {\f{\circdot{x}}{\tau}} = \f{F}{{\tau, {\f{x}{\tau}}}, {\f{u}{\tau}}, s}
    \label{eq:ode}
\end{align}
where $\f{x}{\cdot} \in \R^{n_{x}}$ is the state vector, $\f{u}{\cdot} \in \R^{n_{u}}$ is the control input vector, $s \in \R$ is the parameter, which in our case, is the time-of-flight, and $\f{F}{\cdot} \in \R \times \R^{n_{x}} \times \R^{n_{u}} \times \R \rightarrow \R^{n_{x}}$ is the nonlinear function representing the time-dilated dynamics, which is assumed to be at least once continuously differentiable.

The first-order Taylor series expansion of Equation \eqref{eq:ode} about an arbitrary reference trajectory $\!\left({\f{\overline{x}}{\tau}}, {\f{\overline{u}}{\tau}}, \overline{s}\right)\!$ yields a linear time-varying (LTV) system, given by Equation \eqref{eq:dyn_lin}:
\begin{align}
    \circdot{x}(\tau) \approx \f{A}{\tau}\!\f{x}{\tau} + \f{B}{\tau}\!\f{u}{\tau} + \f{S}{\tau}\!s + \f{d}{\tau} \label{eq:dyn_lin}
\end{align}
where
\begin{subequations}
\begin{align}
    {\f{A}{\tau}} &\defeq \nabla_{\!x\,}{\f{F}{\tau, {\f{\overline{x}}{\tau}}, {\f{\overline{u}}{\tau}}, \overline{s}}}\\
    {\f{B}{\tau}} &\defeq \nabla_{\!u\,}{\f{F}{\tau, {\f{\overline{x}}{\tau}}, {\f{\overline{u}}{\tau}}, \overline{s}}}\\
    {\f{S}{\tau}} &\defeq \nabla_{\!\!\:s\:}{\f{F}{\tau, {\f{\overline{x}}{\tau}}, {\f{\overline{u}}{\tau}}, \overline{s}}}\\
    \begin{split}
    {\f{d}{\tau}} &\defeq {\f{F}{\tau, {\f{\overline{x}}{\tau}}, {\f{\overline{u}}{\tau}}, \overline{s}}} \\
    &\qquad - {\f{A}{\tau}}\,{\f{\overline{x}}{\tau}} - {\f{B}{\tau}}\,{\f{\overline{u}}{\tau}} - {\f{S}{\tau}}\,\overline{s}
    \end{split} \label{eq:w_t_cont_time} 
\end{align}
\label{eqs:jacobians}%
\end{subequations}
Equation \eqref{eq:dyn_lin} is approximate as a consequence of linearization via truncation of the higher-order ($\ge 2$) terms in the Taylor series expansion.

We adopt a first-order hold (FOH) parameterization of the control input signal, which, in contrast to pseudospectral discretization methods, possesses the following characteristics: (i) inter-sample satisfaction of the convex control constraints is guaranteed (provided they are satisfied at the discrete temporal nodes); and, (ii) the resulting conic subproblem has a sparsity pattern that is amenable to real-time implementation \citep{malyuta2019discretization, szmuk2020successive}. These properties make FOH attractive for optimal control applications. 

With FOH, the control input variables are only defined at the discrete temporal nodes, and the continuous-time control input signal is obtained by linearly interpolating between the discrete values at successive nodes. Note that the control input signal is restricted to a continuous piecewise-affine function of time, with only a finite number ($N$) of discrete control variables, as shown in Equation \eqref{eq:control_basis_func}:
\begin{gather}
    \f{u}{\tau} = \f{\sigma_{k}^{-}}{\tau} u_{k} + \f{\sigma_{k}^{+}}{\tau} u_{k+1},\enskip\forall \tau \in \left[\tau_{k}, \tau_{k+1}\right)
    \label{eq:control_basis_func}
\end{gather}
where
\begin{flalign*}
    &\f{\sigma_{k}^{-}}{\tau} \defeq \frac{\tau_{k+1} - \tau}{\tau_{k+1} - \tau_{k}},\enskip\f{\sigma_{k}^{+}}{\tau} \defeq \frac{\tau - \tau_{k}}{\tau_{k+1} - \tau_{k}},\enskip k = \range{1}{N\!-\!1}
\end{flalign*}
The LTV dynamics in terms of deviations from the reference can now be written using the piecewise-affine control input parameterization given by Equation \eqref{eq:control_basis_func}. $\Delta \square$ denotes the deviation of a quantity from its reference, i.e., $\Delta \square \defeq \square - \overline{\square}$, and $\f{\Delta \circdot{x}}{\tau} \defeq \f{\circdot{x}}{\tau} - \f{F}{\tau, {\f{\overline{x}}{\tau}}, {\f{\overline{u}}{\tau}}, \overline{s}}$. Henceforth, ``$=$'' is used in lieu of ``$\approx$'' for notational simplicity, with the understanding that the LTV system is a first-order approximation of the original nonlinear system. The LTV dynamics are given by Equation \eqref{eq:dyn_lin_interp}:
\begin{align}
    \f{\Delta\circdot{x}}{\tau} &= \f{A}{\tau} \Delta \f{x}{\tau} + \f{B}{\tau} \f{\sigma_{k}^{-}}{\tau} \Delta u_{k} + \f{B}{\tau} \f{\sigma_{k}^{+}}{\tau} \Delta u_{k+1} + \f{S}{\tau}\Delta s
    \label{eq:dyn_lin_interp}
\end{align}
Equation \eqref{eq:dyn_lin_interp} has a unique solution \citep{antsaklis2006linear, SCPToolboxCSM2022}, given by Equation \eqref{eq:lin_sys_theory}: $\forall \tau \in \left[\tau_{k}, \tau_{k+1}\right)$, 
\begin{align}
        \Delta\f{x}{\tau} = \f{\Phi}{\tau, \tau_{k}}\f{\Delta x}{\tau_{k}} + \int_{\tau_{k}}^{\tau} \f{\Phi}{\tau, \zeta} \left\{\f{B}{\zeta} \f{\sigma_{k}^{-}}{\zeta} \Delta u_{k} + \f{B}{\zeta} \f{\sigma_{k}^{+}}{\zeta} \Delta u_{k+1} + \f{S}{\zeta} \Delta s \right\} \mathrm{d}\zeta \label{eq:lin_sys_theory}%
\end{align}
where $\f{\Phi}{\tau, \tau_{k}}$, the \textit{state transition matrix}, satisfies the matrix differential equation given by Equation \eqref{eq:mat_diff_eq}:
\begin{gather}
    \f{\circdot{\Phi}}{\tau, \tau_{k}} = \f{A}{\tau}\f{\Phi}{\tau, \tau_{k}}, \enskip \f{\Phi}{\tau_{k}, \tau_{k}} = I_{n_{x}}
    \label{eq:mat_diff_eq}%
\end{gather}

\subsection{Discretization}

\begin{figure}[H]
\centering
    \begin{mybox2}
    \centering
    \includegraphics[width=0.55\linewidth]{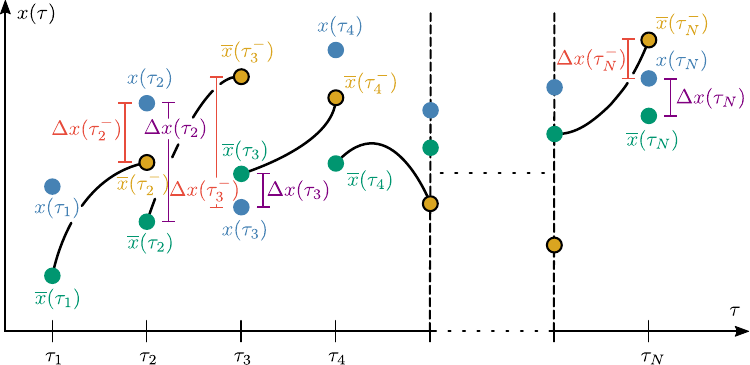}
    \end{mybox2}
    \caption{Propagation of the state in the discretization procedure.}
    \label{fig:stitching}
\end{figure}
Evaluating Equation \eqref{eq:lin_sys_theory} at $\tau = \tau^{-}_{k+1}$, we get Equation \eqref{eq:dyn_disc}:
\begin{align}
\begin{split}
    {\color{bricks}{\f{\Delta x}{\tau^{-}_{k+1}}}} &= A_{k} \f{\Delta x}{\tau_{k}} + B_{k}^{-} \Delta u_{k} + B_{k}^{+} \Delta u_{k+1} + S_{k} \Delta s
\label{eq:dyn_disc}
\end{split}
\end{align}
Note that Equation \eqref{eq:stitching}, which we refer to as the \textit{stitching condition}, holds, as is evident from Figure \ref{fig:stitching}:
\begin{align}
     {\color{bricks}{\f{\Delta x}{\tau^{-}_{k+1}}}} + {\color{goldenrod}{\f{\overline{x}}{\tau^{-}_{k+1}}}} = {\color{violet}{\f{\Delta x}{\tau_{k+1}}}} + {\color{darts}{\f{\overline{x}}{\tau_{k+1}}}} = {\color{steelblue}{\f{x}{\tau_{k+1}}}}, \enskip k = \range{1}{N\!-\!1}\label{eq:stitching}%
\end{align}
The discretized dynamics can now be written as shown in Equation \eqref{eq:dyn_disc_stitch}:
\begin{subequations}
\begin{align}
    {\color{bricks}{\f{\Delta x}{\tau^{-}_{k+1}}}} &= {\color{violet}{\Delta x_{k+1}}} + {\color{darts}{\overline{x}_{k+1}}} - {\color{goldenrod}{x_{k+1}^{\mathrm{prop}}}} \\
    &= A_{k}\,\Delta x_{k} + B_{k}^{-}\,\Delta u_{k} + B_{k}^{+}\,\Delta u_{k+1} + S_{k} \Delta s \\
    \implies \Delta x_{k+1} &= A_{k}\,\Delta x_{k} + B_{k}^{-}\,\Delta u_{k} + B_{k}^{+}\,\Delta u_{k+1} + S_{k}\,\Delta s + d_{k}
\end{align}
\label{eq:dyn_disc_stitch}%
\end{subequations}
where $\Delta x_{k} \defeq \f{\Delta x}{\tau_{k}}$, ${\f{\Delta x}{\tau_{1}}}\!\defeq 0$, ${\color{violet}{\Delta x_{k+1}}} \defeq {\color{violet}{\f{\Delta x}{\tau_{k+1}}}}$, ${\color{goldenrod}{x_{k+1}^{\mathrm{prop}}}} \defeq {\color{goldenrod}{\f{\overline{x}}{\tau^{-}_{k+1}}}}$, ${\color{darts}{\overline{x}_{k+1}}} \defeq {\color{darts}{\f{\overline{x}}{\tau_{k+1}}}}$, and $\f{\overline{u}}{\tau} \defeq \f{\sigma_{k}^{-}}{\tau} \overline{u}_{k} + \f{\sigma_{k}^{+}}{\tau} \overline{u}_{k+1},\,\forall \tau \in \left[\tau_{k}, \tau_{k+1}\right)$.

Here,
\begin{subequations}
\allowdisplaybreaks
\begin{alignat}{1}
    A_{k} &= I_{n_{x}} + \!\!\!\lim_{\:~y \to \tau_{k+1}^{-}}\!\int_{\tau_{k}}^{y} {\f{A}{\zeta}}\,{\f{\Psi_{A}}{\zeta}}\,\mathrm{d}\zeta \\
    B_{k}^{-} &= \!\!\!\lim_{\:~y \to \tau_{k+1}^{-}}\!\int_{\tau_{k}}^{y} \{{\f{A}{\zeta}}\,{\f{\Psi_{B^{-}}}{\zeta}} + {\f{B}{\zeta}}\,\sigma^{-}_{k}(\zeta)\}\,\mathrm{d}\zeta \label{eq:Bkm}\\
    B_{k}^{+} &= \!\!\!\lim_{\:~y \to \tau_{k+1}^{-}}\!\int_{\tau_{k}}^{y}  \{{\f{A}{\zeta}}\,{\f{\Psi_{B^{+}}}{\zeta}} + {\f{B}{\zeta}}\,\sigma^{+}_{k}(\zeta)\}\,\mathrm{d}\zeta \label{eq:Bkp}\\
    S_{k} &= \!\!\!\lim_{\:~y \to \tau_{k+1}^{-}}\!\int_{\tau_{k}}^{y} \{{\f{A}{\zeta}}\,{\f{\Psi_{S}}{\zeta}} + {\f{S}{\zeta}}\}\,\mathrm{d}\zeta \\
    d_{k} &= x_{k+1}^{\mathrm{prop}} - \overline{x}_{k+1}
\end{alignat}
\label{eqs:prop_integrals}
\end{subequations}
where
\begin{subequations}
\begin{align}
    {\f{\circdot{\Psi}_{A}}{\tau}} &= {\f{A}{\tau}}\,{\f{\Psi_{A}}{\tau}} \\
    {\f{\circdot{\Psi}_{B^{-}}}{\tau}} &= {\f{A}{\tau}}\,{\f{\Psi_{B^{-}}}{\tau}} + {\f{B}{\tau}}\,\sigma^{-}_{k}(\tau)\\
    {\f{\circdot{\Psi}_{B^{+}}}{\tau}} &= {\f{A}{\tau}}\,{\f{\Psi_{B^{+}}}{\tau}} + {\f{B}{\tau}}\,\sigma^{+}_{k}(\tau) \\
    {\f{\circdot{\Psi}_{S}}{\tau}} &= {\f{A}{\tau}}\,{\f{\Psi_{S}}{\tau}} + {\f{S}{\tau}}
\end{align}
\label{eqs:disc_ivp}
\end{subequations}
\begin{alignat}{1}
\intertext{\(x_{k+1}^{\mathrm{prop}}\) is evaluated as follows:}
x_{k+1}^{\mathrm{prop}} &= \overline{x}_{k} + \!\!\!\lim_{\:~y \to \tau_{k+1}^{-}}\!\int_{\tau_{k}}^{y} \f{F}{{\zeta, {\f{\overline{x}}{\zeta}}}, {\f{\overline{u}}{\zeta}}, \overline{s}}\,\mathrm{d}\zeta
\label{eqs:prop}
\end{alignat}%
and
${\f{\Psi_{A}}{\tau}} \defeq {\f{\Phi}{\tau, \tau_{k}}}$.


\section{Sequential Conic Optimization (SeCO)}\label{sec:seco}
\vspace{-1em}
\begin{figure}[H]
    \centering
    \includegraphics[width=0.95\linewidth]{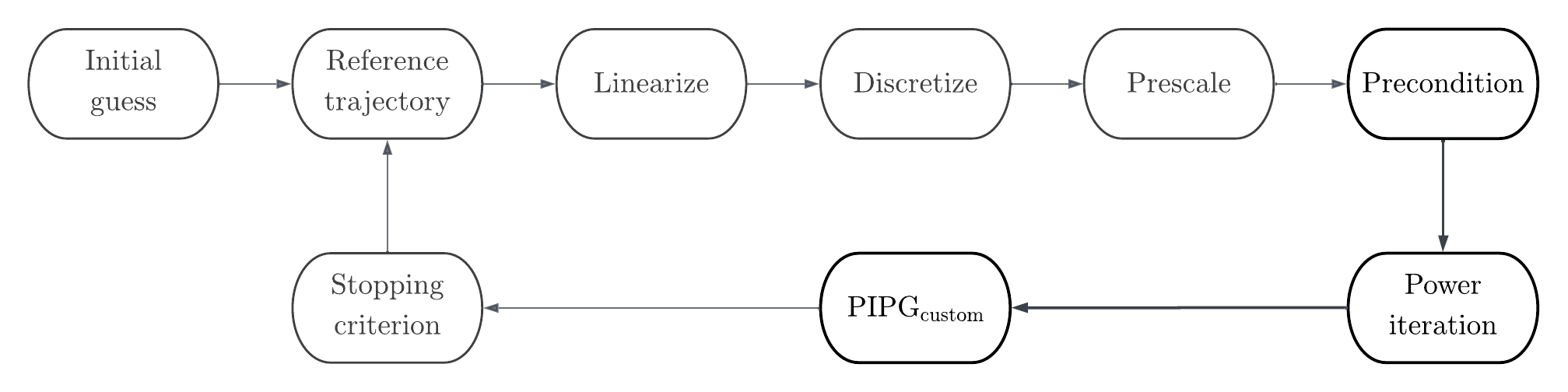}
    \caption{An overview of the SeCO framework; the blocks in bold constitute the low-level solver.}
    \label{fig:SeCO_BD}%
\end{figure}
\vspace{-1em}
In this section, we discuss the components of the sequential conic optimization (\seco{}) framework, shown in Figure \ref{fig:SeCO_BD}, and set up the discrete-time conic subproblem.
\subsection{Virtual State}
\textit{Artificial infeasibility} refers to the phenomenon wherein linearization of the nonconvex constraints in a problem can render the convex subproblem infeasible even if there exists a feasible solution to the original nonconvex problem. In order to mitigate this, one approach is to add slack variables to the linearized constraints, thus ensuring that the subproblem is always feasible. These slack variables are usually unconstrained, yet heavily penalized.

 The value of these slack variables should go to zero at convergence for a solution to be feasible with respect to the original nonconvex problem. The slack variable is usually referred to as \textit{virtual control} when it is added to the linearized dynamics, and \textit{virtual buffer} when it is added to the linearized constraints \citep{SCPToolboxCSM2022}. It has been observed that virtual buffer terms are usually not required if virtual control is used \citep{szmuk2020successive, reynolds2020dual}. In such cases, however, the intermediate reference solutions will not be feasible with respect to the LTV dynamics unless the value of the virtual control is zero. If dynamic feasibility of the intermediate reference trajectories is of importance, one approach could be to buffer the constraints and leave the dynamics equation unchanged.

We opt to use a third approach, called \textit{virtual state} \citep{kamath2023real}. The general idea of this approach is to entirely decouple the dynamics and the control constraints from the state constraints, and to exactly satisfy all the constraints at each iteration, while ensuring that the subproblem is always feasible. In order to achieve this, we introduce a new variable, the virtual state, which acts as a \textit{copy} of the original state variable. If $x$ is the actual state vector, $u$ is the control input vector, and $\xi$ is the virtual state vector, the dynamics constraint is imposed on $x$, the control constraints are imposed on $u$, and all the state constraints are imposed on $\xi$.

In order to ensure that both the dynamics and the path constraints are satisfied at convergence, we minimize the error between $x$ and $\xi$ by including the squared distance between them as a quadratic penalty term in the objective function and heavily penalizing it. The virtual state approach has the benefit of not altering the dynamics manifold, unlike the virtual control approach \citep{szmuk2020successive, reynolds2020dual}, and preserving the shapes of the conic state constraint sets, unlike the virtual buffer approach \citep{SCPToolboxCSM2022}. The left superscript, ${^\xi}\square$, is used to denote virtual state variables.

\subsection{Trust Region}
The trust region radius is the distance between the solution to a subproblem and the trajectory about which the system and nonconvex constraints are linearized to create that subproblem. The purpose of a trust region is to: (1) make sure that the solution does not venture too far from the reference trajectory so as to ensure that the linearization is sufficiently accurate, thus preserving its validity; and, (2) mitigate \textit{artificial unboundedness}, which refers to the phenomenon wherein linearization can render the cost unbounded from below even if it is bounded in the original nonconvex problem \citep{SCPToolboxCSM2022}. 

There exist both hard and soft trust region methods in the literature \citep{szmuk2020successive, mao2016successive}. We impose a soft trust region by augmenting the objective function with a quadratic penalty term, and use the penalized trust region (\ptr) algorithm, which, as the name suggests, penalizes the trust region radius in lieu of constraining it and adopting an outer-loop update rule \citep{szmuk2020successive, reynolds2020real, reynolds2020dual}. This approach has been shown to work very well in practice, and has been employed to successfully solve a wide range of challenging nonconvex optimal control problems in the context of real-time quadrotor path-planning \citep{szmuk2017convexification, szmuk2018real, szmuk2019real}, spacecraft rendezvous and docking \citep{malyuta2020fast, malyuta2021fast, malyuta2021convex}, hypersonic entry trajectory optimization \citep{mceowen2022hypersonic, mceowen2023high}, and real-time rocket landing guidance \cite{reynolds2020dual, reynolds2020real, szmuk2020successive, kim2022guided, kamath2023real}.

\subsection{Initial Guess Generation}

We generate an initial guess trajectory by performing a dual quaternion screw linear interpolation (\sclerp) \cite{kenwright2012dual} between the initial dual quaternion and the nominal final
dual quaternion. A straight-line interpolation is adopted for the mass profile, and the thrust magnitude profile is derived to counteract the (time-varying) weight of the vehicle—which is then saturated based on the thrust magnitude bounds. The gimbal deflection angle, the gimbal azimuth angle, the angular body rates, and the body torques are set to zeros. The linear velocity profile is obtained by linearly interpolating between the initial and terminal conditions. The time-of-flight parameter can be guessed using an analytical procedure such as the one provided in \citep{d1997optimal}. The initial guess generation strategy is a design choice and is highly problem dependent; other approaches can be adopted here as well \cite{malyuta2019discretization, reynolds2020dual, szmuk2020successive}.

\subsection{Prescaling}

The decision variables are prescaled to ensure that the solutions are roughly on the same order of magnitude (between 0 and 1, in our case). This is an important aspect of numerical optimization algorithms, and can significantly improve solution quality and speed up convergence. The interested reader is referred to \citep{SCPToolboxCSM2022} for more details on variable scaling.

\subsection{Constraint Classification} \label{subsec:constraint_classification}

The template {\seco} subproblem is strongly convex, and is given by Equations \eqref{eq:conic_vec}:
\begin{subequations}
\begin{align}
    \underset{z}{\text{minimize}} \quad &\frac{1}{2} z^{\top} Q z + \langle q, z \rangle\\
    \text{subject to} \quad &\begin{aligned}[t]
        & H z - h = 0\\
        & z \in \D
    \end{aligned}
\end{align}
\label{eq:conic_vec}%
\end{subequations}
where all the decision variables are stacked into a single high-dimensional vector, $z$. Here, $Q$ is a positive definite matrix, and $\D$ is a closed convex set that admits a closed-form projection operation. This vectorized problem possesses a sparsity structure that is amenable to customization, as described in Algorithm \ref{alg:pipg_custom}.
 
We prefer to include as many constraints in the form \(z\in\D\) as possible for the following reasons: (1) The constraints \(z\in\D\) will be satisfied at every iteration of {\pipg}. Hence, including more constraints in set \(\D\) leads to a smaller search space. In contrast, the constraint \(Hz-h=0\) is only satisfied asymptotically and does not affect the search space; (2) Transforming many popular constraint sets in optimal control--such as \(\ell_2\)-norm balls, cylinders, boxes--into conic constraints, requires extra auxiliary variables. Imposing these constraints in the form \(z\in\D\) eliminates the need for such auxiliary variables and, as a result, decreases the problem size. As shown in the discrete-time conic subproblem in Section \ref{subsec:discretized_subproblem}, every constraint other than the dynamics is classified into set \(\D\)—this means that these path constraints are exactly satisfied at every solver iteration (up to constraint approximations), and may prove beneficial in the case of premature termination of the solver in an emergency scenario, for instance. The dynamics constraint belongs to the zero cone, and is satisfied asymptotically as the solver converges to an optimum.

\subsection{Constraint Reformulations}

We choose to \textit{combine} certain intersecting path constraints for the following reasons: (1) to reduce the number of constraints imposed (and in turn, reduce the number of operations that the solver needs to carry out); and, (2) to ensure that all the path constraint sets possess closed-form projection operations. Closed-form expressions only exist for the projection onto the intersection of convex sets in special cases---such as the intersection of a cone and a ball, and the intersection of two halfspaces---but this is not the case in general, even if the individual constraint sets can be projected onto \citep{bauschke2018projecting}. Although iterative methods such as the alternating direction method of multipliers (\textsc{admm}) \cite{boyd2011distributed} can be used to compute projections onto the intersection of such convex sets, we opt to reformulate the intersecting constraint sets so as to enable closed-form projections instead.
\subsubsection{Combined Thrust Vector Constraints}
Given the FOH parameterization of the control input, the thrust vector magnitude constraints, in discrete-time, can be expressed as shown in Equations \eqref{eq:control_mag_disc}:
\begin{subequations}
\begin{alignat}{2}
    T_{\min} &\le T_{k} \le T_{\max},\quad &&k = \range{1}{N}\\
    0 &\le \delta_{k} \le \delta_{\max}, &&k = \range{1}{N}\\
    0 &\le \phi_{k} \le 2\pi, &&k = \range{1}{N}
\end{alignat}
\label{eq:control_mag_disc}%
\end{subequations}
where $N$ is the size of the chosen temporal grid, i.e., the number of discrete temporal nodes in the discrete-time subproblem. Further, the thrust vector rate constraints can be expressed as shown in Equations \eqref{eq:control_rate_disc}:
\begin{subequations}
\begin{alignat}{2}
     -\dot{T}_{\max} &\le \frac{T_{k} - T_{k-1}}{(\frac{s}{N-1})} \le \dot{T}_{\max}, \quad &&k = \range{2}{N} \\
     -\dot{\delta}_{\max} &\le \frac{\delta_{k} - \delta_{k-1}}{(\frac{s}{N-1})} \le \dot{\delta}_{\max}, \quad &&k = \range{2}{N} \\
     -\dot{\phi}_{\max} &\le \frac{\phi_{k} - \phi_{k-1}}{(\frac{s}{N-1})} \le \dot{\phi}_{\max}, \quad &&k = \range{2}{N}
\end{alignat}
\label{eq:control_rate_disc}%
\end{subequations}
which are exact, owing to the FOH parameterization \cite{mceowen2023high}. The variable $s$ is the dilation factor, and $\frac{s}{N-1}$ is the length of each of the uniformly spaced time-intervals.

These thrust vector magnitude and rate constraints can be combined by means of an approximation leveraging the reference solution, as shown in Subsection \ref{subsec:discretized_subproblem}, such that $T_{k}$, $\delta_{k}$, and $\phi_{k}$, $k = \range{1}{N}$, are the sole decision variables. This form of the constraints enables closed-form projections, ensures that the original thrust vector constraints are never violated, and is exact at convergence.
\subsubsection{Combined State Constraints}
We propose a new approach to modeling state-triggered constraints (STCs) \cite{szmuk2019scitech, szmuk2020successive, reynolds2019state, malyuta2021advances} that allows for the combination of the global state constraints and the STCs, thus avoiding intersecting constraint sets on the state variables and, in turn, enabling closed-form projection operations onto the constraint sets. For the dual quaternion variable specifically, we reformulate the constraints so as to enable (closed-form) projections onto the intersection of halfspaces.

With the assumption that the global state bound on any given variable, $\square_{\max}$, is greater than its STC counterpart, $\square_{\textsc{stc}_{\max}}$, i.e., $\square_{\max} > \square_{\textsc{stc}_{\max}}$, we observe that the bounds can be expressed in a single expression as follows:
\begin{align}
 g(\square) \le \max\!\left\{-\psi(t)\,\square_{\max},\,\square_{\textsc{stc}_{\max}}\right\}
 \label{eq:combined_state_constraints}%
\end{align}
where $g(\cdot)$ is the constraint function under consideration. The trigger function, $\psi(t)$, given by Equation \eqref{eq:trig_func}:
\begin{align}
    \psi(t) \defeq \operatorname{sgn}\left\{\left(\rho_{\max} - \norm{2\,\dq{q}^{[5:8]}(t)}_{2}\right) \left(\norm{2\,\dq{q}^{[5:8]}(t)}_{2} - \rho_{\min}\right)\right\}\, \ni\, \psi:[0, t_{f}) \to \{-1,\,0,\,+1\}
    \label{eq:trig_func}
\end{align}
takes the value $-1$ if the vehicle is outside the trigger window, $+1$ if the vehicle is inside the trigger window, and $0$, if the vehicle is at either of the triggers. With this formulation, the RHS of Equation \eqref{eq:combined_state_constraints} can automatically switch between the global bound and the STC bound based on the value that the trigger function assumes.

Further, an approximation, similar to the one made with the thrust vector constraints, is made to the trigger function, and the global state constraints and the STCs are combined, as shown in Subsection \ref{subsec:discretized_subproblem}. Note that the combined state constraints are exact at convergence.
\subsection{Projections}
As shown in the discrete-time conic subproblem in Subsection \ref{subsec:discretized_subproblem}, which is a second-order cone program (SOCP), every single path constraint admits a closed-form projection operation. The constraint sets listed in \dx{green} possess direct closed-form projection operations. The ones listed in \fx{blue} and \goldenrod{yellow} involve closed-form projections onto the intersection of halfspaces; the maximum tilt constraint is linearized to enable that. The interested reader is referred to \citep{bauschke2011convex} for a description of these closed-form projection operations.

\vspace{\baselineskip}
\begin{mybox}
\vspace{-0.75em}
\subsection{The Discretized Conic Subproblem}\label{subsec:discretized_subproblem}
$$
\footnotesize
\begin{aligned}
    &\quad\underset{s,\,u_{[1:N]}}{\text{minimize}}~~~ \underbracket[0.2ex]{-w_{m}\, m_N\vphantom{\frac{1}{2}\,w_{\mathrm{vse}} \sum_{k=1}^{N} \norm{x_{k} - \xi_{k}}_{2}^2}}_{\text{cost term}} + \underbracket[0.2ex]{\frac{1}{2}\!\left(w_{\mathrm{tr}}\sum_{k=1}^{N} \left(\|x_{k} - \overline{x}_{k}\|_{2}^2 + \|u_{k} - \overline{u}_{k}\|_{2}^2 \right) + w_{\mathrm{tr}_s}\|s - \overline{s}\|_{2}^{2}\right)}_{\text{soft trust region term}} + \underbracket[0.2ex]{\frac{1}{2}\,w_{\mathrm{vse}} \sum_{k=1}^{N} \norm{x_{k} - \xi_{k}}_{2}^2}_{\text{virtual state penalty term}}\\[0.75em]
    &\begin{aligned}
     \text{subject to}\quad & \fbox{\text{Dynamics}}\\
     \rx{\text{zero cone}}\quad & x_{k+1} = A_{k}\,x_{k} + B_{k}^{-}\,u_{k} + B_{k}^{+}\,u_{k+1} + S_{k}\,s\ + d_{k} &(\mathbb{K}_{0})~&~~ k = \range{1}{N\!-1}\\
     & \fbox{\text{Combined control constraints}}\\
     \dx{\text{box}}\quad & T_{\min} \le T_{1} \le T_{\max} &(\mathbb{D})~&~~\\
     \dx{\text{box}}\quad & 0 \le \delta_{1} \le \delta_{\max} &(\mathbb{D})~&~~\\
     \dx{\text{box}}\quad & 0 \le \phi_{1} \le 2\pi &(\mathbb{D})~&~~\\
     \dx{\text{box}}\quad & \max\!\left(T_{\min},\,-\dot{T}_{\max}\frac{\overline{s}}{N-1} + \overline{T}_{k-1}\right) \le T_{k} \le \min\!\left(T_{\max},\,\dot{T}_{\max}\frac{\overline{s}}{N-1} + \overline{T}_{k-1}\right) &(\mathbb{D})~&~~ k = \range{2}{N}\\
     \dx{\text{box}}\quad & \max\!\left(0,\, -\dot{\delta}_{\max} \frac{\overline{s}}{N-1} + \overline{\delta}_{k-1} \right) \le \delta_{k} \le \min\!\left(\delta_{\max},\, \dot{\delta}_{\max} \frac{\overline{s}}{N-1} + \overline{\delta}_{k-1} \right) &(\mathbb{D})~&~~ k = \range{2}{N}\\
     \dx{\text{box}}\quad & \max\!\left(0,\, -\dot{\phi}_{\max} \frac{\overline{s}}{N-1} + \overline{\phi}_{k-1} \right) \le \phi_{k} \le \min\!\left(2\pi,\, \dot{\phi}_{\max} \frac{\overline{s}}{N-1} + \overline{\phi}_{k-1} \right) &(\mathbb{D})~&~~ k = \range{2}{N}\\
     \dx{\text{box}}\quad & \norm{\tau_{k}}_{\infty} \leq \tau_{\max} &(\mathbb{D})~&~~ k = \range{1}{N}\\
     & \fbox{\text{Combined state constraints}} & & \\
     \dx{\text{box}}\quad & \norm{^\xi\!\!\>\dq{\omega}^{[1:3]}_{k}}_{\infty} \le \max\!\left(-\overline{\psi}_{k}\,\omega_{\max},\,\omega_{\textsc{stc}_{\max}}\right) & (\mathbb{D})~&~~ k = \range{2}{N\!-1}\\
     \dx{\text{ball}}\quad & \norm{^\xi\!\!\>\dq{\omega}^{[4:6]}_{k}}_2 \le \max\!\left(-\overline{\psi}_{k}\,v_{\max},\,v_{\textsc{stc}_{\max}}\right) & (\mathbb{D})~&~~ k = \range{2}{N\!-1}\\
     \hphantom{\quad}\fx{\text{halfspace}}\quad & \overline{\dq{q}}^{[1:2]^{\top}}_{k} \!^\xi\!\!\>\dq{q}_{k}^{[1:2]} \le \norm{\overline{\dq{q}}^{[1:2]}_{k}}_{2}\sin\!\tfrac{\theta_{\textsc{stc}_{\max}}}{2} &(\mathbb{D})~&~~ k \ni \overline{\psi}_{k} \ge 0\\
     \fx{\text{halfspace}}\quad & \text{Linearize Equation \eqref{eq:LoS_ncvx} (maximum line-of-sight angle)} &(\mathbb{D})~&~~ k \ni \overline{\psi}_{k} \ge 0\\
     \goldenrod{\text{halfspace}}\quad & \overline{\dq{q}}^{[1:2]^{\top}}_{k} \!^\xi\!\!\>\dq{q}_{k}^{[1:2]} \le \norm{\overline{\dq{q}}^{[1:2]}_{k}}_{2}\sin\!\tfrac{\theta_{\max}}{2} &(\mathbb{D})~&~~ k \ni \overline{\psi}_{k} < 0\\
     \goldenrod{\text{halfspace}}\quad & \text{Linearize Equation \eqref{eq:min_alt_ncvx} (minimum altitude)} & (\mathbb{D})~&~~ k \ni \overline{\psi}_{k} < 0\\
     & \fbox{\text{Boundary conditions}}\\
     \dx{\text{singleton}}\quad & ^\xi\!\!\>m_{1} = m_{i} & (\mathbb{D})~& \\
     \dx{\text{singleton}}\quad & ^\xi\!\!\>\dq{q}_{1} = \begin{pmatrix}\q{q}_{i} \\ \dfrac{1}{2}\!\begin{pmatrix}\q{r}_{\I_{i}}\\ 0\end{pmatrix}\otimes\q{q}_{i}\end{pmatrix} & (\mathbb{D})~& \\
     \dx{\text{singleton}}\quad & ^\xi\!\!\>\dq{\omega}_{1} = \begin{pmatrix}\q{\omega}_{\B_{i}}\\q_{i}^{*} \otimes \begin{pmatrix}\q{v}_{\I_{i}}\\ 0\end{pmatrix} \otimes q_{i}\end{pmatrix} & (\mathbb{D})~\\
     \dx{\text{halfspace}}\quad & ^\xi\!\!\>m_{N} \geq m_{f} & (\mathbb{D})~& \\
     \dx{\text{singleton}}\quad & ^\xi\!\!\>\dq{q}^{[1:2]}_{N} = 0 & (\mathbb{D})~& \\
     \dx{\text{subspace}}\quad & \begin{bmatrix} \left(\frac{1}{2}[r_{\mathcal{I}_f}]_{\otimes}\!\begin{pmatrix}\mathbf{0}_{2\times2}\\ I_{2}\end{pmatrix}\right) & -I_{4} \end{bmatrix} {}^\xi\!\!\>\dq{q}^{[3:8]}_{N} = 0 & (\mathbb{D})~& \\
      \dx{\text{singleton}}\quad & ^\xi\!\!\>\dq{\omega}_{N} = \begin{pmatrix}
          \mathbf{0}_{5 \times 1} \\ v_{z_{\I_{f}}}
      \end{pmatrix} & (\mathbb{D})~&\\
     \end{aligned}
\end{aligned}
$$
\vspace{0em}
\end{mybox}
\section{High-Performance Solver}\label{sec:solver}
\vspace{-0.5em}
\begin{figure}[H]
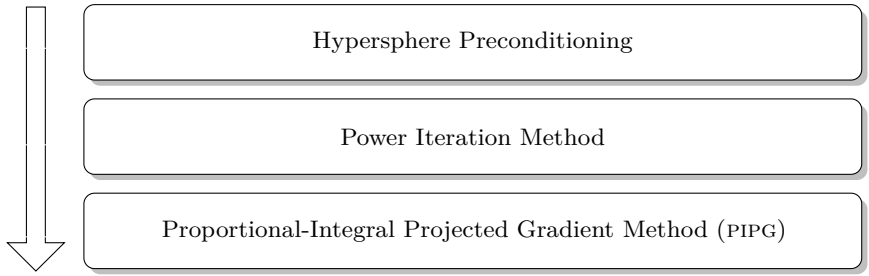

\centering
\begin{minipage}[b]{.98\linewidth}
\tikzset{priority arrow fill/.style={
fill=white}}
\tikzset{priority arrow/.style={
draw=black,
single arrow,
minimum height=\distancemodules,
minimum width=0.75cm,
priority arrow fill,
rotate=90,
single arrow head extend=0.25cm,
anchor=west}}
\tikzset{priority arrow/.append style={rotate=180,anchor=0,xshift=1,}}
\smartdiagramset{border color=black, 
                 set color list={white,white,white},
                 description title font=\scriptsize,
                 descriptive items y sep=1.25cm,
                 description title text width=8cm,
                 description title width=8cm,
                 description width=8cm,
                 description text width=10cm,
                 uniform arrow color=true,
                 arrow color=black,
                 priority arrow width=1.25cm,
                 priority arrow head extend=0.25cm,
                 priority arrow height advance=1cm,
                 priority tick size=0cm}
\centering
\smartdiagram[priority descriptive diagram]{
Proportional-Integral Projected Gradient Method (\pipg),
Power Iteration Method,
Hypersphere Preconditioning}
\caption{The SeCO subproblem solver.}
\label{fig:solver_levels}
\end{minipage}
\end{figure}
\vspace{-0.5em}

In this section, we describe the development of a high-performance first-order conic optimization solver, as outlined in Figure \ref{fig:solver_levels}, to solve the convex subproblems in \seco{}, which we then customize in Section \ref{sec:custom}.
\subsection{Preconditioning}

First-order methods are sensitive to problem conditioning and typically perform poorly on ill-conditioned problems \citep{stellato2020osqp}. As SCP (\seco) implementations typically impose a heavy weight on the virtual term penalty relative to the original cost and trust region penalties \citep{SCPToolboxCSM2022}, the objective function of the resulting convex (conic) subproblem is inherently ill-conditioned.

Preconditioning is a heuristic that seeks to reduce the condition number of a parameter matrix in an optimization problem via a coordinate transformation. There exist both exact and heuristic methods for preconditioning, both with their own limitations. The exact optimal diagonal preconditioner for a matrix, i.e., the diagonal preconditioner that minimizes the condition number of the resulting matrix, can be found by solving a semidefinite program (SDP) \citep{stellato2020osqp, boyd1994linear}. However, such an approach is not suitable for real-time applications for the following reasons: (1) solving such a problem is typically more computationally expensive than solving the original problem itself \citep{stellato2020osqp}; and, (2) SDPs are generally more difficult to solve than second-order cone programs (SOCPs), which usually form the most general class of convex optimization problems that are suitable for safety-critical applications \citep{malyuta2021advances}.

Matrix equilibration is another approach that has been widely used in the literature for preconditioning in the context of optimization \citep{ruiz2001scaling, sinkhorn1967concerning, pock2011diagonal, giselsson2014diagonal, fougner2018parameter}. This approach typically involves finding a diagonal preconditioner to scale the matrix under consideration such that its rows have equal norms and its columns have equal norms \citep{diamond2017stochastic}. Equilibration-based preconditioning techniques either rely on convex optimization or iterative methods, some more reliable than others \citep{stellato2020osqp}. Although many of these methods take both the objective function and constraints into account and have been found to work well in practice \citep{stellato2020osqp, o2016conic}, they either require the solution of another convex optimization problem or rely on heuristic iterative procedures that are not guaranteed to reduce the condition number of the matrix being equilibrated even if they converge. Further, preconditioners based on incomplete matrix factorizations are also popular in practice \cite{trefethen1997numerical}.

\begin{figure}[!htb]
    \centering
    \begin{minipage}[t]{0.4775\linewidth}
        \begin{mybox}
        \centering
        \includegraphics[width=0.75\linewidth]{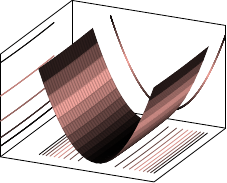}
        \end{mybox}
        \caption{An ill-conditioned bivariate quadratic function (prior to preconditioning, for instance).}
        \label{fig:ill}
    \end{minipage}\hspace*{\fill}
    \begin{minipage}[t]{0.4775\linewidth}
        \begin{mybox}
        \centering
        \includegraphics[width=0.75\linewidth]{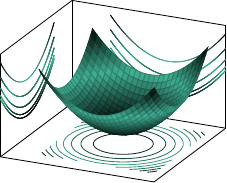}
        \end{mybox}
        \caption{A well-conditioned bivariate quadratic function (after preconditioning, for instance).}
        \label{fig:well}
    \end{minipage}
\end{figure}
In this work, we use the hypersphere preconditioner \cite{kamath2025optimal}, shown in Algorithm \ref{alg:precondition}, for problems that fit the template of {\seco}. This procedure uses an exact Cholesky factorization of the objective function, and a row-normalization of the constraint matrix, while effectively exploiting the structure of the conic subproblem in {\seco} with minimal computational overhead—in fact, the Cholesky factor can be computed analytically, without the need for any explicit matrix factorization/inversion operations or iterative heuristics; i.e., the structure of the {\seco} subproblem enables the computation of the Cholesky factor of the objective function matrix in closed-form, as shown in the Algorithm \ref{alg:precondition_custom}. Further, to compute the optimal cost-scaling factor in terms of minimizing the condition number of the Karush Kuhn Tucker (KKT) matrix of the problem \cite{kamath2025optimal}, we adopt the shifted power iteration method, the customized version of which we describe in Algorithm \ref{alg:shifted_power_custom}. Further, the transformed dynamics matrix still possesses a sparsity structure that is amenable to customization, as shown in Algorithm \ref{alg:pipg_custom}. The preconditioner transforms ill-conditioned quadratic functions to well-conditioned quadratic functions, the level sets of which are $\ell_{2}$-norm hyperspheres; hence the name. The effect of preconditioning on an ill-conditioned quadratic objective function, such as the one in Figure \ref{fig:ill}, is depicted in Figure \ref{fig:well}. In our initial testing, we observed a reduction in the number of solver iterations to convergence of roughly $5$ to $10$ times with this preconditioner. Further, we normalize the rows of the constraint matrix—this helps reduce the condition number of the preconditioned constraint matrix in practice \citep{kamath2025optimal}.
\begin{algorithm}[H]
\small
\caption{Hypersphere Preconditioning}\label{alg:precondition}
    \vspace{0.25em}
    \begin{flushleft}
        \textbf{Inputs:} $Q$, $q$, $H$, $h$, $\D$
    \end{flushleft}
    \begin{algorithmic}[1]
    \Require $Q \succ 0$, $Q = Q^{\top}$
    \vspace{1em}
    \State $L^{\top} L \leftarrow \chol Q$ \Comment{Cholesky decomposition}
    \State $L_{\mathrm{inv}} \leftarrow L^{-1}$
    \State $\hat{\D} \leftarrow L\,\D$
    \State $\hat{H} \leftarrow H L_{\mathrm{inv}}$
    \State $\hat{H} \leftarrow E \hat{H}$ \Comment{$E$: diagonal matrix with reciprocals of row-norms of $\hat{H}$}
    \State $\hat{h} \leftarrow E h$
    \State $\sigma_{\max} \leftarrow$ Algorithm \ref{alg:power} \Comment{power iteration}
    \State $\sigma_{\min} \leftarrow$ Algorithm \ref{alg:shifted_power} \Comment{shifted power iteration}
    \State $\lambda \leftarrow \sqrt{\frac{\sigma_{\min}}{2}}$
    \State $\hat{q} \leftarrow \lambda\,L_{\mathrm{inv}}^{\top} q$\vspace{1em}
    \Ensure $\square \in \hat{\D} \Leftrightarrow L^{-1}\,\square \in \D$
    \end{algorithmic}
    \begin{flushleft}
        \textbf{Return:} $\lambda$, $\hat{q}$, $\hat{H}$, $\hat{h}$, $\hat{\D}$, $L$, $L_{\mathrm{inv}}$, $\sigma_{\max}$
    \end{flushleft}
    \vspace{0.25em}
\end{algorithm}
As shown in Algorithm \ref{alg:precondition}, we scale the objective function with a positive scalar, $\lambda$, that factors into the step-sizes of {\pipg} and hence serves as a tuning parameter. Scaling the objective function appropriately can also help keep the magnitude of the dual variables in check \citep{stellato2020osqp}. 
The particular choice of $\lambda$ in Algorithm \ref{alg:precondition} leads to minimization of the condition number of the KKT matrix of the preconditioned problem \cite{kamath2025optimal}. That said, $\lambda$ can be manually tuned to obtain good performance as well.

The preconditioned conic subproblem is shown in Problem \ref{prob:conic_preconditioned}:
\begin{subequations}
\begin{align}
    \underset{\hat{z}}{\text{minimize}} \quad & \frac{\lambda}{2}\,\hat{z}^{\top}\hat{z} + \langle \hat{q}, \hat{z} \rangle\\
    \text{subject to} \quad &\begin{aligned}[t]
        & \hat{H} \hat{z} - \hat{h} = 0\\
        & \hat{z} \in \hat{\D}
    \end{aligned}
\label{prob:conic_preconditioned}
\end{align}
\end{subequations}
The objective function matrix of the preconditioned problem is of the form $\lambda\,I$, with a condition number of $1$, which is the minimum attainable condition number (and hence, the preconditioner is optimal). The transformed projection $\hat{z} \in \hat{\D}$ needs to preserve membership of the original primal variable, $z$, to set $\D$, i.e., $\hat{z} \in \hat{\D} \Leftrightarrow z \in \D$ \citep{o2016conic}. The structure of the conic subproblem in {\seco} naturally preserves this, as is apparent from Algorithm \ref{alg:precondition_custom}. Post-solution, the original primal variable can be recovered as follows: $z = L_{\mathrm{inv}} \hat{z}$. 
\subsection{Power Iteration Method}

The power iteration method described in Algorithm \ref{alg:power} computes the maximum singular value of a given diagonalizable matrix 
\cite{trefethen1997numerical}. We can exploit the structure of matrix $\hat{H}$ in trajectory optimization problems to customize Algorithm \ref{alg:power} to Algorithm \ref{alg:power_custom} so as to avoid sparse linear algebra operations with large dimensional matrices and vectors.

\begin{algorithm}[H]
\small
\caption{Power Iteration Method}\label{alg:power}
    \vspace{0.25em}
    \begin{flushleft}
        \textbf{Inputs:} $\hat{H}$, $z$, $\epsilon_{\mathrm{abs}}$, $\epsilon_{\mathrm{rel}}$, $\epsilon_{\mathrm{buff}}$, $j_{\max}$
    \end{flushleft}
    \begin{algorithmic}[1]
    \Require $\norm{z}_{2} > 0$
    \vspace{1em}
    \State $\sigma \leftarrow \norm{z}_{2}$\label{line:init} \Comment{initialization}\vspace{1ex}
    \For {$j \leftarrow \range{1}{j_{\max}}$}\vspace{1ex}
    \State $w \leftarrow \frac{1}{\sigma} \hat{H} z$\label{line:w}
    \State $z \leftarrow \hat{H}^{\top} w$\label{line:z}
    \State $\sigma^{\star} \leftarrow \norm{z}_{2}$\label{line:sigma_star}\vspace{1ex}
    \If {$\abs{\sigma^{\star} - \sigma} \le \epsilon_{\mathrm{abs}} + \epsilon_{\mathrm{rel}}\,\max\{\sigma^{\star},\,\sigma\}$} \Comment{stopping criterion}
    \State \textbf{break}
    \ElsIf {$j < j_{\max}$}
    \State $\sigma \leftarrow \sigma^{\star}$
    \EndIf\vspace{1ex}
    \EndFor\vspace{1ex}
    \State $\sigma \leftarrow (1 + \epsilon_{\mathrm{buff}})\,\sigma^{\star}$ \Comment{buffer the (under) estimated maximum singular value}\vspace{0.125em}
    \end{algorithmic}
    \begin{flushleft}
        \textbf{Return:} $\sigma$ \Comment{$\approx \max \operatorname{spec} \hat{H}^{\top}\hat{H} = \sigma_{\max}(\hat{H}^{\top}\hat{H}) = \|\hat{H}\|_{2}^{2}$}
    \end{flushleft}
    \vspace{0.25em}
\end{algorithm}

\subsection{Shifted Power Iteration Method}

The shifted power iteration method described in Algorithm \ref{alg:shifted_power} is similar to the power iteration method, except that it computes an estimate for the smallest singular value of a given diagonalizable matrix instead \cite{wilkinson1988algebraic}. See the appendix in \cite{kamath2025optimal} for a detailed description of this algorithm. Note that the shifted power iteration method is used solely to estimate the smallest singular value of $\hat{H}\,\hat{H}^{\top}$ to determine the optimal value for the cost-scaling factor, $\lambda$. If $\lambda$ is manually tuned instead, this algorithm is not required. The customized version of this algorithm is given by Algorithm \ref{alg:shifted_power_custom}.

\begin{algorithm}[H]
\small
\caption{Shifted Power Iteration Method}\label{alg:shifted_power}
    \vspace{0.5em}
    \begin{flushleft}
        \textbf{Inputs:} $\hat{H}$, $w$, $\epsilon_{\mathrm{abs}}$, $\epsilon_{\mathrm{rel}}$, $\epsilon_{\mathrm{buff}}$, $j_{\max}$, $\sigma_{\max}$
    \end{flushleft}
    \begin{algorithmic}[1]
    \Require $\norm{w}_{2} > 0$
    \vspace{1em}
    \State $\tilde{\sigma} \leftarrow \norm{w}_{2}$\label{line:init_shifted} \Comment{initialization}\vspace{1ex}
    \For {$j \leftarrow \range{1}{j_{\max}}$}\vspace{1ex}
    \State $z \leftarrow \hat{H}^{\top}w$\label{line:z_shifted}\vspace{1ex}
    \State $w \leftarrow \frac{1}{\tilde{\sigma}}\,(\hat{H}\,z - \sigma_{\max}\,w)$\label{line:w_shifted} \Comment{$\sigma_{\max} \defeq \|\hat{H}\|^{2}$}\vspace{1ex}
    \State $\tilde{\sigma}^{\star} \leftarrow \norm{w}_{2}$\label{line:sigma_star_shifted}\vspace{1ex}
    \If {$\abs{\tilde{\sigma}^{\star} - \tilde{\sigma}} \le \epsilon_{\mathrm{abs}} + \epsilon_{\mathrm{rel}}\,\max\{\tilde{\sigma}^{\star},\,\tilde{\sigma}\}$} \Comment{stopping criterion}
    \State \textbf{break}
    \ElsIf {$j < j_{\max}$}
    \State $\tilde{\sigma} \leftarrow \tilde{\sigma}^{\star}$
    \EndIf\vspace{1ex}
    \EndFor\vspace{1ex}
    \State $\sigma_{\min} \leftarrow (1 - \epsilon_{\mathrm{buff}})\,(\sigma_{\max} - \tilde{\sigma}^{\star})$
    \Comment{buffer the (over) estimated minimum singular value}\vspace{0.125em}
    \end{algorithmic}
    \begin{flushleft}
        \textbf{Return:} $\sigma_{\min}$ \Comment{$\approx \min \operatorname{spec} \hat{H}\,\hat{H}^{\top} = \sigma_{\min}(\hat{H}\,\hat{H}^{\top})$}
    \end{flushleft}
\end{algorithm}

\newpage
\subsection{PIPG}

The proportional-integral projected gradient method, {\pipg}, is a first-order primal-dual algorithm for conic optimization \citep{yu2021proportional}. It allows matrix-factorization/inverse-free and easily-verifiable solver implementations for real-time and embedded applications \cite{yu2022real}. {\pipg} achieves the optimal global convergence rates (worst-case) in theory, and performs much faster in practice \citep{yu2022extrapolated}. It exploits the sparsity structure of conic constraints via parallel matrix operations and the geometric structure of constraint sets via efficient closed-form projections \cite{yu2020proportional}. Unlike most off-the-shelf methods, {\pipg} allows for warm-starting and enjoys a light computational overhead, as it avoids the cumbersome canonical transformation procedure that standard conic programs are subject to. It is also compatible with extrapolation, which has been shown to accelerate convergence \citep{yu2022extrapolated}. {\pipg} not only enables versatile convex optimization, but also has the ability to boost the performance of sequential convex programming (SCP) methods for nonconvex optimization. Moreover, the {\seco} framework specializes SCP algorithms to exploit features of {\pipg} to solve nonconvex optimal control problems in real-time \cite{kamath2023real}.

One of the major applications of {\pipg} is in solving trajectory optimization problems, given that all of the sparse linear algebra operations in Equations \ref{alg:pipg} can be devectorized \cite[Algorithm 2]{elango2022customized}, as shown in Algorithm \ref{alg:pipg_custom}, and posed as simple, small-dimensional matrix-vector manipulations that are suitable for real-time performance onboard resource-constrained embedded hardware.

{\pipg} converges to an optimal solution when the difference between two consecutive iterates converges to zero \citep[Theorem 1]{yu2022extrapolated}. Hence, in practice, we terminate the solver when the difference between two consecutive iterates is sufficiently small. In particular, given a relative accuracy tolerance \(\epsilon_{\text{rel}}\) and an absolute accuracy tolerance \(\epsilon_{\text{abs}}\), we terminate {\pipg} when the following conditions are met: 
\begin{align*}
\begin{split}
    \norm{\hat{z}^{j+1} - \hat{z}^{j}}_{\infty} &\le \epsilon_{\text{abs}} + \epsilon_{\text{rel}}\,\max\left\{\norm{\hat{z}^{j+1}}_{\infty},\,\norm{\hat{z}^{j}}_{\infty}\right\},\\
    \norm{\hat{w}^{j+1} - \hat{w}^{j}}_{\infty} &\le \epsilon_{\text{abs}} + \epsilon_{\text{rel}}\,\max\left\{\norm{\hat{w}^{j+1}}_{\infty},\,\norm{\hat{w}^{j}}_{\infty}\right\}
\end{split}    
\end{align*}
We note that such a combination of absolute and relative accuracy tolerances is popular among first-order solvers \citep{stellato2020osqp,o2021operator}.


\begin{algorithm}[H]
\small
\caption{\pipg}\label{alg:pipg}
    \vspace{0.25em}
    \begin{flushleft}
        \textbf{Inputs:} $\hat{q}$, $\hat{H}$, $\hat{h}$, $\hat{\D}$, $L$, $L_{\mathrm{inv}}$, $\lambda$, $\sigma$, $\rho$, $\epsilon_{\mathrm{abs}}$, $\epsilon_{\mathrm{rel}}$, $j_{\mathrm{check}}$, $j_{\max}$,\\[1ex]
        \hphantom{\textbf{Inputs:}\,} $z^{\star}$, $\hat{w}^{\star}$ \Comment{warm start}
    \end{flushleft}
    \vspace{0.125em}
    \begin{algorithmic}[1]
    \State $\hat{z}^{\star} \leftarrow L\,z^{\star}$ \Comment{transform previous primal solution}
    \State $\zeta^{1} \leftarrow \hat{z}^{\star}$ \Comment{initialize transformed primal variable}
    \State $\eta^{1} \leftarrow \hat{w}^{\star}$ \Comment{initialize transformed dual variable}\vspace{1ex}
    \State $\alpha \leftarrow \frac{2}{\lambda + \sqrt{\lambda^{2} + 4\sigma}}$ \Comment{step-size}\vspace{1ex}
    \For {$j \leftarrow \range{1}{j_{\max}}$}\vspace{1ex}
    \State $\hat{z}^{j+1} = \pi_{\D}[\zeta^{j}-\alpha\,(\lambda\,\zeta^{j} + \hat{q} + \hat{H}^{\top} \eta^{j})]$ \Comment{projected gradient step}
    \State $\hat{w}^{j+1} = \eta^{j} + \alpha\,(\hat{H}(2\,\hat{z}^{j+1}-\zeta^{j}) - \hat{h})$ \Comment{PI feedback of affine equality constraint violation}
    \State $\zeta^{j+1} = (1 - \rho)\,\zeta^{j} + \rho\,\hat{z}^{j+1}$ \Comment{extrapolate transformed primal variable}
    \State $\eta^{j+1} = (1 - \rho)\,\eta^{j} + \rho\,\hat{w}^{j+1}$ \Comment{extrapolate transformed dual variable}\vspace{1ex}
    \If {$j \operatorname{mod} j_{\mathrm{check}} = 0$} \Comment{check stopping criterion every $j_{\mathrm{check}}$ iterations}\vspace{1ex}
    \If {$\norm{\hat{z}^{j+1} - \hat{z}^{j}}_{\infty} \le \epsilon_{\text{abs}} + \epsilon_{\text{rel}}\,\max\!\left\{\norm{\hat{z}^{j+1}}_{\infty},\,\norm{\hat{z}^{j}}_{\infty}\right\}$ \textbf{and}\\
    $\qquad\quad\enskip\, \norm{\hat{w}^{j+1} - \hat{w}^{j}}_{\infty} \le \epsilon_{\text{abs}} + \epsilon_{\text{rel}}\,\max\!\left\{\norm{\hat{w}^{j+1}}_{\infty},\,\norm{\hat{w}^{j}}_{\infty}\right\}$} \Comment{stopping criterion}\vspace{1ex}
    \State \textbf{break}\vspace{1ex}
    \EndIf\vspace{1ex}
    \EndIf\vspace{1ex}
    \EndFor\vspace{1ex}
    \State $z^{\star} \leftarrow L_{\mathrm{inv}}\,\hat{z}^{j+1}$ \Comment{recover original primal variable}
    \State $\hat{w}^{\star} \leftarrow \hat{w}^{j+1}$ \Comment{retain transformed dual variable}\vspace{0.125em}
    \end{algorithmic}
    \begin{flushleft}
        \textbf{Return:} $z^{\star}$, $\hat{w}^{\star}$
    \end{flushleft}
    \vspace{0.25em}
\end{algorithm}
\section{Solver Customization}\label{sec:custom}
\vspace{-0.5em}
\begin{figure}[H]
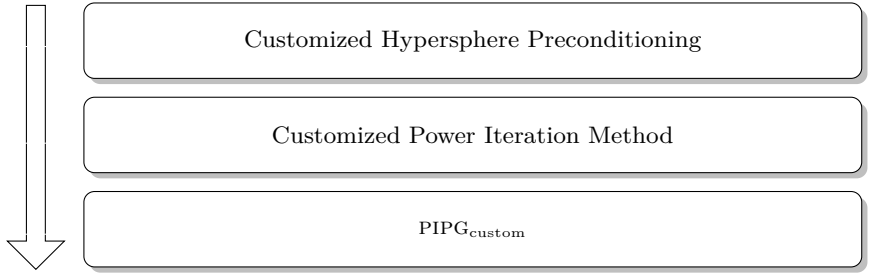

\centering
\begin{minipage}[b]{.98\linewidth}
\tikzset{priority arrow fill/.style={
fill=white}}
\tikzset{priority arrow/.style={
draw=black,
single arrow,
minimum height=\distancemodules,
minimum width=0.75cm,
priority arrow fill,
rotate=90,
single arrow head extend=0.25cm,
anchor=west}}
\tikzset{priority arrow/.append style={rotate=180,anchor=0,xshift=1,}}
\smartdiagramset{border color=black, 
                 set color list={white,white,white},
                 description title font=\scriptsize,
                 descriptive items y sep=1.25cm,
                 description title text width=8cm,
                 description title width=8cm,
                 description width=8cm,
                 description text width=10cm,
                 uniform arrow color=true,
                 arrow color=black,
                 priority arrow width=1.25cm,
                 priority arrow head extend=0.25cm,
                 priority arrow height advance=1cm,
                 priority tick size=0cm}
\centering
\smartdiagram[priority descriptive diagram]{
{\pipgc},
Customized Power Iteration Method,
Customized Hypersphere Preconditioning}
\caption{The customized SeCO subproblem solver.}
\label{fig:custom_solver_levels}
\end{minipage}
\end{figure}
\vspace{-0.5em}

We define customization as the exploitation of the sparsity pattern of the optimal control problem at hand, so as to enable low-dimensional matrix-vector multiplications and other dense linear algebra operations with devectorized variables, and thus avoid sparse linear algebra operations. In contrast to dense matrix-vector multiplications, sparse matrix-vector multiplications (SpMV) typically suffer from: (i) additional computational overhead in terms of instructions and storage; (ii) memory access patterns that are indirect and irregular; and, (iii) more cache misses, and are hence inefficient \citep{williams2007optimization, chari2025qoco}. Customized algorithms preclude these inefficient operations, and are especially effective with optimization problem sizes that are characteristic of onboard guidance \citep{dueri2017customized, elango2022customized, chari2025qoco}. In this section, we describe customization of the algorithms presented in Section \ref{sec:solver}, leading to the customized subproblem solver outlined in Figure \ref{fig:custom_solver_levels}.
\subsection{Customized Preconditioning}
Consider the vectorized conic optimization problem given by Equations \eqref{eq:conic_vec}, where
\begin{gather*}
z = \left(x_{1}, \ldots,\,x_{N},\,\xi_{1}, \ldots,\,\xi_{N},\,u_{1}, \ldots,\,u_{N},\,s\right) \in \R^{n_{z}}\\
q = \left(q_{x},\,q_{\xi},\,q_{u},\,q_{s}\right) \in \R^{n_{z}}\\
Q = \begin{pmatrix}Q_{\text{state}} && \\ & Q_{u} & \\ && Q_{s}\end{pmatrix} \in \mathbb{S}_{++}^{n_{z}}
\end{gather*}
$Q_{\text{state}} \defeq W_{\text{state}} \otimes I_{n_{x} N}$, $Q_{u} \defeq w_{\mathrm{tr}}\,I_{n_{u} N}$, $Q_{s} \defeq w_{\mathrm{tr}_{s}}$, $W_{\text{state}} \defeq \begin{pmatrix}w_{\mathrm{tr}} + w_{\mathrm{vse}} & -w_{\mathrm{vse}} \\ -w_{\mathrm{vse}} & w_{\mathrm{vse}} \end{pmatrix}$, and $n_{z}$ is the length of $z$.

$Q$ is symmetric positive definite (SPD) and therefore has a unique Cholesky decomposition \citep[Corollary 7.2.9]{horn2012matrix}. Further, since $Q$ is block diagonal, the blocks are effectively decoupled, and the Cholesky decomposition can be applied directly to the individual blocks, as shown in Equation \eqref{eq:chol_individual}:
\begin{gather}
    \chol Q = \begin{pmatrix}\chol Q_{\text{state}} && \\ & \chol Q_{u} & \\ && \chol Q_{s}\end{pmatrix}
    \label{eq:chol_individual}
\end{gather}

Let $L$ denote the Cholesky factor of $Q$, such that $Q = L^{\top} L$, $L_{\square}$ denote the Cholesky factor of a block partition of $Q$, such that $Q_{\square} = L_{\square}^{\top} L_{\square}$, and $R$ denote the Cholesky factor of $W_{\text{state}}$, such that $W_{\text{state}} = R^{\top} R$.

Given two SPD matrices $A$ and $B$, it can be shown that $\chol(A \otimes B) = \chol A \otimes \chol B$ \citep{schacke2004kronecker}. Therefore,
\begin{align}
    \chol Q_{\text{state}} &= \chol W_{\text{state}} \otimes \chol I_{n_{x} N}\\
    &= R^{\top} R \otimes I_{n_{x} N}\\
    &= L_{\text{state}}^{\top} L_{\text{state}}
\end{align}
where
\begin{gather}
    R = \frac{1}{\sqrt{w_{\mathrm{tr}} + w_{\mathrm{vse}}}}\begin{pmatrix}w_{\mathrm{tr}} + w_{\mathrm{vse}} & -w_{\mathrm{vse}} \\ 0 & \sqrt{w_{\mathrm{tr}} w_{\mathrm{vse}}}\end{pmatrix} \in \R^{2 \times 2} \label{eq:R}\\
    \text{and}\enskip L_{\text{state}} = R \otimes I_{n_{x} N}
\end{gather}
Since $Q_{u}$ is a diagonal matrix, $\chol Q_{u} = L_{u}^{\top} L_{u}$, where $L_{u} = \sqrt{w_{\mathrm{tr}}}\,I_{n_{u} N}$. Since $Q_{s}$ is a scalar, $\chol Q_{s} = L_{s}^{2}$, where $L_{s} = \sqrt{w_{\mathrm{tr}_{s}}}$. 

Finally, the Cholesky factor of $Q$ can be computed as follows: $L = \blkdiag\{L_{\text{state}}, L_{u}, L_{s}\}$. Since $L$ is block diagonal, its inverse can be written in terms of the inverses of the individual block partitions $L_{\text{state}}$, $L_{u}$, and $L_{s}$, as shown in Equation \eqref{eq:L_inv}:
\begin{gather}
    L^{-1} = \blkdiag\{L_{\text{state}}^{-1}, L_{u}^{-1}, L_{s}^{-1}\}
    \label{eq:L_inv}
\end{gather}

Given two nonsingular SPD matrices $A$ and $B$, it can be shown that $\left(A \otimes B\right)^{-1} = A^{-1} \otimes B^{-1}$ \citep[Corollary 10]{broxson2006kronecker}. Therefore,
\begin{gather}
    L_{\text{state}}^{-1} = \left(R \otimes I_{n_{x} N}\right)^{-1} = R^{-1} \otimes I_{n_{x} N}
\end{gather}
where
\begin{gather}
    R^{-1} = \frac{1}{\sqrt{w_{\mathrm{tr}} + w_{\mathrm{vse}}}}\begin{pmatrix}1 & w \\ 0 & w + \frac{1}{w} \end{pmatrix} \in \R^{2 \times 2} \label{eq:R_inv}\\
    w = \sqrt{\frac{w_{\mathrm{vse}}}{w_{\mathrm{tr}}}}
\end{gather}
and
\begin{align}
    L_{u}^{-1} = \frac{1}{\sqrt{w_{\mathrm{tr}}}} I_{n_{u} N}~\text{and}~L_{s}^{-1} = \frac{1}{\sqrt{w_{\mathrm{tr}_{s}}}}
\end{align}
For embedded applications, if the problem data does not need to change, i.e., if $Q$ is fixed, the preconditioning parameters can be computed offline and stored onboard. However, even in cases where the problem data may change, note that effectively, the only matrix factorization/inversion operations required in the proposed preconditioning procedure are one Cholesky decomposition of a $2 \times 2$ matrix ($W_{\text{state}}$) and one inversion of a $2 \times 2$ upper triangular matrix ($R$), both of which have closed-form expressions. An efficient implementation of the customized hypersphere preconditioning procedure with no explicit matrix factorizations/inversions is documented in Algorithm \ref{alg:precondition_custom}. Note that all the transformations in the algorithm only involve scaling the problem data by scalars (with the exception of one vector addition operation). Further, these scaling factors only depend on the objective function weights, which are independent of the problem size, thus making the algorithm suitable for large-scale problems as well.

Further, given the optimal control structure, row normalization of the constraint matrix can also be customized, as shown in Algorithm \ref{alg:precondition_custom}.

\begin{algorithm}[H]
\small
\caption{Customized Hypersphere Preconditioning}\label{alg:precondition_custom}
    \vspace{0.25em}
    \begin{flushleft}
        \textbf{Inputs:} $w_{\mathrm{vse}}$, $w_{\mathrm{tr}}$, $w_{\mathrm{tr}_{s}}$, $q_{x}$, $q_{\xi}$, $q_{u}$, $q_{s}$, $A^{-}_{[1:N-1]}$, $E^{-}_{[1:N-1]}$, $B^{-}_{[1:N-1]}$, $B^{+}_{[1:N-1]}$, $S_{[1:N-1]}$, $d_{[1:N-1]}$, $\D_{\xi}$, $\D_{u}$, $\D_{s}$, $\lambda$
    \end{flushleft}
    \begin{algorithmic}[1]
    \Require $w_\mathrm{vse},\, w_{\mathrm{tr}},\, w_{\mathrm{tr}_{s}} > 0$ \hfill{$R$: Equation \eqref{eq:R}; $R^{-1}$: Equation \eqref{eq:R_inv}}
    \vspace{1em}
    \State $l_{x_{1}} \leftarrow \sqrt{w_{\mathrm{tr}} + w_{\mathrm{vse}}}$ \Comment{$R_{\{1, 1\}}$}
    \State $l_{x_{2}} \leftarrow \frac{-w_{\mathrm{vse}}}{l_{x_{1}}}$ \Comment{$R_{\{1, 2\}}$}
    \State $l_{\xi} \leftarrow \frac{\sqrt{w_{\mathrm{tr}}\,w_{\mathrm{vse}}}}{l_{x_{1}}}$ \Comment{$R_{\{2, 2\}}$}
    \State $l_{x_{{1}_{\mathrm{inv}}}} \leftarrow \frac{1}{l_{x_{1}}}$ \Comment{$R^{-1}_{\{1, 1\}} = \frac{1}{R_{\{1, 1\}}}$}
    \State $l_{x_{{2}_{\mathrm{inv}}}} \leftarrow \frac{-l_{x_{2}}}{l_{x_{1}} l_{\xi}}$ \Comment{$R^{-1}_{\{1, 2\}} = \frac{-R_{\{1, 2\}}}{R_{\{1, 1\}} R_{\{2, 2\}}}$}
    \State $l_{\xi_{\mathrm{inv}}} \leftarrow \frac{1}{l_{\xi}}$ \Comment{$R^{-1}_{\{2, 2\}} = \frac{1}{R_{\{2, 2\}}}$}
    \State $l_{u} \leftarrow \sqrt{w_{\mathrm{tr}}}$
    \State $l_{u_\mathrm{inv}} \leftarrow \frac{1}{l_{u}}$
    \State $l_{s} \leftarrow \sqrt{w_{\mathrm{tr}_{s}}}$
    \State $l_{s_\mathrm{inv}} \leftarrow \frac{1}{l_{s}}$
    \State $\hat{\D}_{\xi} \leftarrow l_{\xi}\,\D_{\xi}$
    \State $\hat{\D}_{u} \leftarrow l_{u}\,\D_{u}$
    \State $\hat{\D}_{s} \leftarrow l_{s}\,\D_{s}$
    \State $\hat{A}^{-}_{[1:N-1]} \leftarrow l_{x_{{1}_{\mathrm{inv}}}} A_{[1:N-1]}$
    \State $\hat{A}^{+}_{[1:N-1]} \leftarrow -l_{x_{{1}_{\mathrm{inv}}}} I_{n_{x}}$
    \State $\hat{E}^{-}_{[1:N-1]} \leftarrow l_{x_{{2}_{\mathrm{inv}}}} A_{[1:N-1]}$
    \State $\hat{E}^{+}_{[1:N-1]} \leftarrow -l_{x_{{2}_{\mathrm{inv}}}} I_{n_{x}}$
    \State $\hat{B}^{-}_{[1:N-1]} \leftarrow l_{u_{\mathrm{inv}}} B^{-}_{[1:N-1]}$
    \State $\hat{B}^{+}_{[1:N-1]} \leftarrow l_{u_{\mathrm{inv}}} B^{+}_{[1:N-1]}$
    \State $\hat{S}_{[1:N-1]} \leftarrow l_{s_{\mathrm{inv}}} S_{[1:N-1]}$
    \For{$k = \range{1}{N-1}$} \Comment{row normalization}
        \For{$l = \range{1}{n_{x}}$}
            \State $r \leftarrow \left[\hat{A}^{-}_{[k]}[l, :], \hat{A}^{+}_{[k]}[l, :], \hat{E}^{-}_{[k]}[l, :], \hat{E}^{+}_{[k]}[l, :], \hat{B}^{-}_{[k]}[l, :], \hat{B}^{+}_{[k]}[l, :], \hat{S}_{[k]}[l, :]\right]$
            \State $n \leftarrow \norm{r}_{\infty}$
            \State $r \leftarrow \frac{1}{n} r$
            \State $\hat{d}[l, :] \leftarrow \frac{1}{n} d[l, :]$
        \EndFor
    \EndFor
    \State $\sigma_{\max} \leftarrow$ Algorithm \ref{alg:power_custom}
    \State $\sigma_{\min} \leftarrow$ Algorithm \ref{alg:shifted_power_custom}
    \State $\lambda \leftarrow \sqrt{\frac{\sigma_{\min}}{2}}$
    \State $\hat{q}_{x} \leftarrow \lambda\,l_{x_{{1}_{\mathrm{inv}}}} q_{x}$
    \State $\hat{q}_{\xi} \leftarrow \lambda\,(l_{x_{{2}_{\mathrm{inv}}}} q_{x} + l_{\xi_\mathrm{inv}} q_{\xi})$
    \State $\hat{q}_{u} \leftarrow \lambda\,l_{u_\mathrm{inv}} q_{u}$
    \State $\hat{q}_{s} \leftarrow \lambda\,l_{s_\mathrm{inv}} q_{s}$\vspace{0.125em}
    \end{algorithmic}
    \begin{flushleft}
        \textbf{Return:} $\hat{q}_{x}$, $\hat{q}_{\xi}$, $\hat{q}_{u}$, $\hat{q}_{s}$, $\hat{A}^{-}_{[1:N-1]}$, $\hat{A}^{+}_{[1:N-1]}$, $\hat{E}^{-}_{[1:N-1]}$, $\hat{E}^{+}_{[1:N-1]}$, $\hat{B}^{-}_{[1:N-1]}$, $\hat{B}^{+}_{[1:N-1]}$, $\hat{S}_{[1:N-1]}$, $\hat{d}_{[1:N-1]}$,\\
        \hphantom{\textbf{Return:}} $\hat{\D}_{\xi}$, $\hat{\D}_{u}$, $\hat{\D}_{s}$, $l_{x_{1}}$, $l_{x_{2}}$, $l_{\xi}$, $l_{u}$, $l_{s}$, $l_{x_{{1}_{\mathrm{inv}}}}$, $l_{x_{{2}_{\mathrm{inv}}}}$, $l_{\xi_{\mathrm{inv}}}$, $l_{u_{\mathrm{inv}}}$, $l_{s_{\mathrm{inv}}}$, $\sigma_{\max}$
    \end{flushleft}
    \vspace{0.25em}
\end{algorithm}

\subsection{Customized Power Iteration Method}

\begin{algorithm}[H]
\small
\caption{Customized Power Iteration Method}
\label{alg:power_custom}
    \vspace{0.25em}
    \begin{flushleft}
        \textbf{Inputs:} $\hat{A}^{-}_{[1:N-1]}$, $\hat{A}^{+}_{[1:N-1]}$, $\hat{E}^{-}_{[1:N-1]}$, $\hat{E}^{+}_{[1:N-1]}$, $\hat{B}^{-}_{[1:N-1]}$, $\hat{B}^{+}_{[1:N-1]}$, $\hat{S}_{[1:N-1]}$,\\[1ex]
        \hphantom{\textbf{Inputs:}\,} $\epsilon_{\mathrm{abs}}$, $\epsilon_{\mathrm{rel}}$, $\epsilon_{\mathrm{buff}}$, $j_{\max}$,\\[1ex]
        \hphantom{\textbf{Inputs:}\,} $x_{[1:N]}$, $\xi_{[1:N]}$, $u_{[1:N]}$, $s$, $w_{[1:N-1]}$
    \end{flushleft}
    \begin{algorithmic}[1]
    \Require $\norm{x_{[1:N]}}_{2} > 0$,\, $\norm{\xi_{[1:N]}}_{2} > 0$,\, $\norm{u_{[1:N]}}_{2} > 0$,\, $s > 0$
    \vspace{1em}
    \State $\sigma \leftarrow s^{2}$\vspace{1ex}
    \For {$k \leftarrow \range{1}{N}$}\Comment{Algorithm \ref{alg:power}, Line \ref{line:init}}
    \State $\sigma \leftarrow \sigma + \norm{x_{k}}_{2}^{2} + \norm{\xi_{k}}_{2}^{2} + \norm{u_{k}}_{2}^{2}$
    \EndFor\vspace{1ex}
    \State $\sigma \leftarrow \sqrt{\sigma}$\vspace{1ex}
    \For {$j \leftarrow \range{1}{j_{\max}}$}\vspace{1ex}
    \For {$k \leftarrow \range{1}{N\!-\!1}$}\Comment{Algorithm \ref{alg:power}, Line \ref{line:w}}\vspace{1ex}
    \State $w_{k} \leftarrow \frac{1}{\sigma}\!\left(\hat{A}^{-}_{k}\,x_{k} + \hat{A}^{+}_{k}\,x_{k+1} + \hat{E}^{-}_{k}\,\xi_{k} + \hat{E}^{+}_{k}\,\xi_{k+1} + \hat{B}^{-}_{k} u_{k} + \hat{B}^{+}_{k} u_{k+1} + \hat{S}_{k}\,s\right)$\vspace{1ex}
     \EndFor\vspace{1ex}
    \State $x_{1} \leftarrow \hat{A}^{-^{\top}}_{1}\!w_{1}$
    \State $\xi_{1} \leftarrow \hat{E}^{-^{\top}}_{1}\!w_{1}$
    \State $u_{1} \leftarrow \hat{B}^{-^{\top}}_{1}\!w_{1}$
    \State $s \leftarrow \hat{S}^{\top}_{1}w_{1}$\vspace{1ex}
    \For {$k \leftarrow \range{2}{N\!-\!1}$}\Comment{Algorithm \ref{alg:power}, Line \ref{line:z}}
    \State $x_{k} \leftarrow \hat{A}^{-^{\top}}_{k}\!w_{k} + \hat{A}^{+^{\top}}_{k-1}\!w_{k-1}$
    \State $\xi_{k} \leftarrow \hat{E}^{-^{\top}}_{k}\!w_{k} + \hat{E}^{+^{\top}}_{k-1}\!w_{k-1}$
    \State $u_{k} \leftarrow \hat{B}^{-^{\top}}_{k}\!w_{k} + \hat{B}^{+^{\top}}_{k-1}w_{k-1}$
    \State $s \leftarrow s + \hat{S}^{\top}_{k}w_{k}$\vspace{0.5ex}
    \EndFor\vspace{1ex}
    \State $x_{N} \leftarrow \hat{A}^{+^{\top}}_{N-1}w_{N-1}$
    \State $\xi_{N} \leftarrow \hat{E}^{+^{\top}}_{N-1}w_{N-1}$
    \State $u_{N} \leftarrow \hat{B}^{+^{\top}}_{N-1}w_{N-1}$\vspace{1ex}
    \State $\sigma^{\star} \leftarrow s^{2}$\vspace{1ex}
    \For {$k \leftarrow \range{1}{N}$}\Comment{Algorithm \ref{alg:power}, Line \ref{line:sigma_star}}
    \State $\sigma^{\star} \leftarrow \sigma^{\star} + \norm{x_{k}}_{2}^{2} + \norm{\xi_{k}}_{2}^{2} + \norm{u_{k}}_{2}^{2}$
    \EndFor\vspace{1ex}
    \State $\sigma^{\star} \leftarrow \sqrt{\sigma^{\star}}$\vspace{1ex}
    \If {$\abs{\sigma^{\star} - \sigma} \le \epsilon_{\mathrm{abs}} + \epsilon_{\mathrm{rel}}\,\max\{\sigma^{\star},\,\sigma\}$} \Comment{stopping criterion}
    \State \textbf{break}
    \ElsIf {$j < j_{\max}$}
    \State $\sigma \leftarrow \sigma^{\star}$
    \EndIf\vspace{1ex}
    \EndFor\vspace{1ex}
    \State $\sigma \leftarrow (1 + \epsilon_{\mathrm{buff}})\,\sigma^{\star}$\Comment{buffer the (under) estimated maximum singular value}\vspace{0.125em}
    \end{algorithmic}
    \begin{flushleft}
        \textbf{Return:} $\sigma$\Comment{$\approx \max \operatorname{spec} \hat{H}^{\top}\hat{H} = \sigma_{\max}(\hat{H}^{\top}\hat{H}) = \|\hat{H}\|_{2}^{2}$}
    \end{flushleft}
    \vspace{0.25em}
\end{algorithm}

\subsection{Customized Shifted Power Iteration Method}

\begin{algorithm}[H]
\small
\caption{Customized Shifted Power Iteration Method}
\label{alg:shifted_power_custom}
    \vspace{0.25em}
    \begin{flushleft}
        \textbf{Inputs:} $\hat{A}^{-}_{[1:N-1]}$, $\hat{A}^{+}_{[1:N-1]}$, $\hat{E}^{-}_{[1:N-1]}$, $\hat{E}^{+}_{[1:N-1]}$, $\hat{B}^{-}_{[1:N-1]}$, $\hat{B}^{+}_{[1:N-1]}$, $\hat{S}_{[1:N-1]}$,\\[1ex]
        \hphantom{\textbf{Inputs:}\,} $\epsilon_{\mathrm{abs}}$, $\epsilon_{\mathrm{rel}}$, $\epsilon_{\mathrm{buff}}$, $j_{\max}$, $\sigma_{\max}$\\[1ex]
        \hphantom{\textbf{Inputs:}\,} $x_{[1:N]}$, $\xi_{[1:N]}$, $u_{[1:N]}$, $s$, $w_{[1:N-1]}$
    \end{flushleft}
    \begin{algorithmic}[1]
    \Require $\norm{w_{[1:N-1]}}_{2} > 0$
    \vspace{1em}
    \State $\tilde{\sigma} \leftarrow \norm{w_{[1:N-1]}}_{2}$\vspace{1ex}
    \For {$j \leftarrow \range{1}{j_{\max}}$}\vspace{1ex}
    \State $x_{1} \leftarrow \hat{A}^{-^{\top}}_{1}\!w_{1}$
    \State $\xi_{1} \leftarrow \hat{E}^{-^{\top}}_{1}\!w_{1}$
    \State $u_{1} \leftarrow \hat{B}^{-^{\top}}_{1}\!w_{1}$
    \State $s \leftarrow \hat{S}^{\top}_{1}w_{1}$\vspace{1ex}
    \For {$k \leftarrow \range{2}{N\!-\!1}$}\Comment{Algorithm \ref{alg:shifted_power}, Line \ref{line:z_shifted}}
    \State $x_{k} \leftarrow \hat{A}^{-^{\top}}_{k}\!w_{k} + \hat{A}^{+^{\top}}_{k-1}\!w_{k-1}$
    \State $\xi_{k} \leftarrow \hat{E}^{-^{\top}}_{k}\!w_{k} + \hat{E}^{+^{\top}}_{k-1}\!w_{k-1}$
    \State $u_{k} \leftarrow \hat{B}^{-^{\top}}_{k}\!w_{k} + \hat{B}^{+^{\top}}_{k-1}w_{k-1}$
    \State $s \leftarrow s + \hat{S}^{\top}_{k}w_{k}$\vspace{0.5ex}
    \EndFor\vspace{1ex}
    \State $x_{N} \leftarrow \hat{A}^{+^{\top}}_{N-1}w_{N-1}$
    \State $\xi_{N} \leftarrow \hat{E}^{+^{\top}}_{N-1}w_{N-1}$
    \State $u_{N} \leftarrow \hat{B}^{+^{\top}}_{N-1}w_{N-1}$\vspace{1ex}
    \For {$k \leftarrow \range{1}{N\!-\!1}$}\Comment{Algorithm \ref{alg:shifted_power}, Line \ref{line:w_shifted}}\vspace{1ex}
    \State $w_{k} \leftarrow \frac{1}{\tilde{\sigma}}\!\left(\hat{A}^{-}_{k}\,x_{k} + \hat{A}^{+}_{k}\,x_{k+1} + \hat{E}^{-}_{k}\,\xi_{k} + \hat{E}^{+}_{k}\,\xi_{k+1} + \hat{B}^{-}_{k} u_{k} + \hat{B}^{+}_{k} u_{k+1} + \hat{S}_{k}\,s - \sigma_{\max}\,w_{k}\right)$\vspace{1ex}
     \EndFor\vspace{1ex}
    \State $\tilde{\sigma}^{\star} \leftarrow \norm{w_{[1:N-1]}}_{2}$\vspace{1ex}\Comment{Algorithm \ref{alg:shifted_power}, Line \ref{line:sigma_star_shifted}}
    \If {$\abs{\tilde{\sigma}^{\star} - \tilde{\sigma}} \le \epsilon_{\mathrm{abs}} + \epsilon_{\mathrm{rel}}\,\max\{\tilde{\sigma}^{\star},\,\tilde{\sigma}\}$} \Comment{stopping criterion}
    \State \textbf{break}
    \ElsIf {$j < j_{\max}$}
    \State $\tilde{\sigma} \leftarrow \tilde{\sigma}^{\star}$
    \EndIf\vspace{1ex}
    \EndFor\vspace{1ex}
    \State $\sigma_{\min} \leftarrow (1 - \epsilon_{\mathrm{buff}})\,(\sigma_{\max} - \tilde{\sigma}^{\star})$\Comment{buffer the (over) estimated minimum singular value}\vspace{0.125em}
    \end{algorithmic}
    \begin{flushleft}
        \textbf{Return:} $\sigma_{\min}$\Comment{$\approx \min \operatorname{spec} \hat{H}\,\hat{H}^{\top} = \sigma_{\min}(\hat{H}\,\hat{H}^{\top})$}
    \end{flushleft}
    \vspace{0.25em}
\end{algorithm}

\newgeometry{left=0.75in,
             right=0.75in,
             top=0.4in,
             bottom=0.5in,
             footskip=0in}
\subsection{Customized PIPG}
\vspace{-0.5em}
\begin{algorithm}[H]
\footnotesize
\caption{\pipgc{}}\label{alg:pipg_custom}
    \vspace{1ex}
    \begin{flushleft}
        \textbf{Inputs:} $q_{x_{[1:N]}}$, $q_{\xi_{[1:N]}}$, $q_{u_{[1:N]}}$, $q_{s}$,\\ 
        \hphantom{\textbf{Inputs:}} $\hat{A}^{-}_{[1:N-1]}$, $\hat{A}^{+}_{[1:N-1]}$, $\hat{E}^{-}_{[1:N-1]}$, $\hat{E}^{+}_{[1:N-1]}$, $\hat{B}^{-}_{[1:N-1]}$, $\hat{B}^{+}_{[1:N-1]}$, $\hat{d}_{[2:N]}$, $\hat{\overline{x}}_{[1:N]}$, $\hat{\overline{u}}_{[1:N]}$, $\hat{\overline{s}}$,\\
        \hphantom{\textbf{Inputs:}} $\hat{\D}_{x_{1}}$, $\hat{\D}_{\xi_{[2:N]}}$, $\hat{\D}_{u_{[1:N]}}$, $\hat{\D}_{s}$,\\
        \hphantom{\textbf{Inputs:}} $l_{x_{1}}$, $l_{x_{2}}$, $l_{\xi}$, $l_{u}$, $l_{s}$, $l_{x_{1_{\mathrm{inv}}}}$, $l_{x_{2_{\mathrm{inv}}}}$, $l_{\xi_{\mathrm{inv}}}$, $l_{u_{\mathrm{inv}}}$, $l_{s_{\mathrm{inv}}}$,\\
        \hphantom{\textbf{Inputs:}} $\lambda$, $\sigma$, $\omega$, $\rho$, $\epsilon_{\mathrm{abs}}$, $\epsilon_{\mathrm{rel}}$, $j_{\mathrm{check}}$, $j_{\max}$,\\
        \hphantom{\textbf{Inputs:}} $\Delta\hat{x}_{[1:N]}^{\star}$, $\Delta\hat{\xi}_{[1:N]}^{\star}$, $\Delta\hat{u}_{[1:N]}^{\star}$, $\Delta\hat{s}^{\star}$, $w^{\star}_{[1:N-1]}$\Comment{warm start}
    \end{flushleft}
    \vspace{-0.75em}
    \begin{algorithmic}[1]
    \State $\Delta x_{\zeta_{[1:N]}}^{1} \leftarrow l_{x_{1}}\,\Delta \hat{x}_{[1:N]}^{\star} + l_{x_{2}}\,\Delta\hat{\xi}_{[1:N]}^{\star}$\Comment{initialize primal variables}
    \State $\Delta \xi_{\zeta_{[1:N]}}^{1} \leftarrow l_{\xi}\,\Delta \hat{\xi}_{[1:N]}^{\star}$
    \State $\Delta u_{\zeta_{[1:N]}}^{1} \leftarrow l_{u}\,\Delta \hat{u}_{[1:N]}^{\star}$
    \State $\Delta s_{\zeta}^{1} \leftarrow l_{s}\,\Delta \hat{s}^{\star}$
    \State $\eta_{[1:N-1]}^{1} \leftarrow w_{[1:N-1]}^{\star}$\Comment{initialize dual variable}\vspace{1ex}
    \State $\alpha \leftarrow \frac{2}{\lambda + \sqrt{\lambda^{2} + 4\omega\sigma}}$\Comment{step-sizes} 
    \State $\beta \leftarrow \omega\alpha$\vspace{1ex}
    \For {$j \leftarrow \range{1}{j_{\max}}$}\vspace{1ex}
    \State $\Delta \hat{x}^{j+1}_{1} \leftarrow \pi_{\hat{\D}_{x_{1}}}[\Delta x_{\zeta_{1}}^{j}-\alpha\,(\lambda\,\Delta x_{\zeta_{1}}^{j} + \lambda\,q_{x_{1}} + \hat{A}^{-^{\top}}_{1}\!\eta_{1}^{j}) + \hat{\overline{x}}^{j}_{1}] - \hat{\overline{x}}^{j}_{1}$
    \State $\Delta \hat{\xi}^{j+1}_{1} \leftarrow 0$
    \State $\Delta \hat{u}^{j+1}_{1} \leftarrow \pi_{\hat{\D}_{u_{1}}}[\Delta u_{\zeta_{1}}^{j}-\alpha\,(\lambda\,\Delta u_{\zeta_{1}}^{j} + \lambda\,q_{u_{1}} + \hat{B}^{-^{\top}}_{1}\!\eta^{j}_{1}) + \hat{\overline{u}}^{j}_{1}] - \hat{\overline{u}}^{j}_{1}$
    \State $\Delta \mathcal{S} \leftarrow \hat{S}^{\top}_{1}\eta^{j}_{1}$\vspace{1ex}
    \For {$k \leftarrow \range{2}{N\!-\!1}$}\Comment{projected gradient step}\vspace{1ex}
    \State $\Delta \hat{x}_{k}^{j+1} \leftarrow \Delta x_{\zeta_{k}}^{j}-\alpha\,(\lambda\,\Delta x_{\zeta_{k}}^{j} + \lambda\,q_{x_{k}} + \hat{A}^{-^{\top}}_{k}\!\eta_{k}^{j} + \hat{A}^{+^{\top}}_{k-1}\eta_{k-1}^{j})$
    \State $\Delta \hat{\xi}_{k}^{j+1} \leftarrow \pi_{\hat{\D}_{\xi_{k}}}[\Delta \xi_{\zeta_{k}}^{j}-\alpha\,(\lambda\,\Delta \xi_{\zeta_{k}}^{j} + \lambda\,q_{\xi_{k}} + \hat{E}^{-^{\top}}_{k}\!\eta_{k}^{j} + \hat{E}^{+^{\top}}_{k-1}\eta_{k-1}^{j}) + \hat{\overline{x}}^{j}_{k}] - \hat{\overline{x}}^{j}_{k}$
    \State $\Delta \hat{u}_{k}^{j+1} \leftarrow \pi_{\hat{\D}_{u_{k}}}[\Delta u_{\zeta_{k}}^{j}-\alpha\,(\lambda\,\Delta u_{\zeta_{k}}^{j} + \lambda\,q_{u_{k}} + \hat{B}^{-^{\top}}_{k}\!\eta_{k}^{j} + \hat{B}^{+^{\top}}_{k-1}\eta_{k-1}^{j}) + \hat{\overline{u}}_{k}^{j}] - \hat{\overline{u}}_{k}^{j}$
    \State $\Delta \mathcal{S} \leftarrow \Delta \mathcal{S} + \hat{S}^{\top}_{k}\eta_{k}^{j}$\vspace{1ex}
    \EndFor\vspace{1ex}
    \State $\Delta \hat{x}_{N}^{j+1} \leftarrow \Delta x_{\zeta_{N}}^{j}-\alpha\,(\lambda\,\Delta x_{\zeta_{N}}^{j} + \lambda\,q_{x_{N}} + \hat{A}^{+^{\top}}_{N-1}\eta_{N-1}^{j})$\vspace{1ex}
    \State $\Delta \hat{\xi}_{N}^{j+1} \leftarrow \pi_{\hat{\D}_{\xi_{N}}}[\Delta \xi_{\zeta_{N}}^{j}-\alpha\,(\lambda\,\Delta \xi_{\zeta_{N}}^{j} + \lambda\,q_{\xi_{N}} + \hat{E}^{+^{\top}}_{N-1}\eta_{N-1}^{j}) + \hat{\overline{x}}^{j}_{N}] - \hat{\overline{x}}^{j}_{N}$
    \State $\Delta \hat{u}_{N}^{j+1} \leftarrow \pi_{\hat{\D}_{u_{N}}}[\Delta u_{\zeta_{N}}^{j}-\alpha\,(\lambda\,\Delta u_{\zeta_{N}}^{j} + \lambda\,q_{u_{N}} + \hat{B}^{+^{\top}}_{N-1}\eta_{N-1}^{j}) + \hat{\overline{u}}_{N}^{j}] - \hat{\overline{u}}_{N}^{j}$\vspace{1ex}
    \State $\Delta \hat{s}^{j+1} \leftarrow \pi_{\hat{\D}_{s}}[\Delta s_{\zeta}^{j}-\alpha\,(\lambda\,\Delta s_{\zeta}^{j} + \lambda\,q_{s} + \Delta\mathcal{S}) + \hat{\overline{s}}^{j}] - \hat{\overline{s}}^{j}$\vspace{1ex}
    \For {$k \leftarrow 1:N\!-\!1$}\Comment{PI feedback of affine equality constraint violation}\vspace{1ex}
    \State $w_{k}^{j+1} \leftarrow \eta_{k}^{j} + \beta\,(\hat{A}^{-}_{k}\,(2 \Delta\hat{x}^{j+1}_{k} - \Delta x^{j}_{\zeta_{k}}) + \hat{A}^{+}_{k}\,(2 \Delta\hat{x}^{j+1}_{k+1} - \Delta x^{j}_{\zeta_{k+1}}) + \hat{E}^{-}_{k}\,(2 \Delta\hat{\xi}^{j+1}_{k} - \Delta \xi^{j}_{\zeta_{k}}) + \hat{E}^{+}_{k}\,(2 \Delta\hat{\xi}^{j+1}_{k+1} - \Delta \xi^{j}_{\zeta_{k+1}})$ \\
    \hphantom{$\quad\quad\!\: w_{k}^{j+1} \leftarrow \eta_{k}^{j} + \beta\,($}$+\,\hat{B}^{-}_{k}\,(2 \Delta\hat{u}^{j+1}_{k} - \Delta u^{j}_{\zeta_{k}}) + \hat{B}^{+}_{k}\,(2 \Delta\hat{u}^{j+1}_{k+1} - \Delta u^{j}_{\zeta_{k+1}}) + \hat{S}_{k}\,(2 \Delta\hat{s}^{j+1} - \Delta s^{j}_{\zeta}) + \hat{d}_{k+1})$\vspace{1ex}
    \EndFor\vspace{1ex}
    \State $\Delta x_{\zeta_{[1:N]}}^{j+1} \leftarrow (1 - \rho)\,\Delta x_{\zeta_{[1:N]}}^{j} + \rho\,\Delta \hat{x}_{[1:N]}^{j+1}$ \Comment{extrapolate primal variables}
    \State $\Delta \xi_{\zeta_{[1:N]}}^{j+1} \leftarrow (1 - \rho)\,\Delta \xi_{\zeta_{[1:N]}}^{j} + \rho\,\Delta \hat{\xi}_{[1:N]}^{j+1}$
    \State $\Delta u_{\zeta_{[1:N]}}^{j+1} \leftarrow (1 - \rho)\,\Delta u_{\zeta_{[1:N]}}^{j} + \rho\,\Delta \hat{u}_{[1:N]}^{j+1}$
    \State $\eta_{[1:N-1]}^{j+1} \leftarrow (1 - \rho)\,\eta_{[1:N-1]}^{j} + \rho\,w_{[1:N-1]}^{j+1}$\Comment{extrapolate dual variables}\vspace{1ex}
    \If {$j \operatorname{mod} j_{\mathrm{check}} = 0$} \Comment{check stopping criterion every $j_{\mathrm{check}}$ iterations}\vspace{1ex}
    \State $\textsc{terminate} \leftarrow \textsc{stopping}(\Delta\hat{x}_{[1:N]}^{j+1},\,\Delta\hat{\xi}_{[1:N]}^{j+1},\,\Delta\hat{u}_{[1:N]}^{j+1},\,\Delta\hat{s}^{j+1},\,w_{[1:N-1]}^{j+1}$,
    \State \hphantom{$\textsc{terminate} \leftarrow \textsc{stopping}\:$}$\,\Delta\hat{x}_{[1:N]}^{j},\,\Delta\hat{\xi}_{[1:N]}^{j},\,\Delta\hat{u}_{[1:N]}^{j},\,\,\,\,\Delta\hat{s}^{j\hphantom{+1}}\!\!\!,\,w_{[1:N-1]}^{j},\,\epsilon_{\mathrm{abs}},\,\epsilon_{\mathrm{rel}})$\vspace{1ex}
    \If {$\textsc{terminate} = \textsc{true}$}\vspace{1ex}\Comment{stopping criterion}
    \State \textbf{break}\vspace{1ex}
    \EndIf\vspace{1ex}
    \EndIf\vspace{1ex}
    \EndFor\vspace{1ex}
    \State $\Delta \hat{x}_{[1:N]}^{\star} \leftarrow l_{x_{1_{\mathrm{inv}}}}\,\Delta \hat{x}_{[1:N]}^{j+1} + l_{x_{2_{\mathrm{inv}}}}\,\Delta \hat{\xi}_{[1:N]}^{j+1}$\Comment{update primal variables}
    \State $\Delta \hat{\xi}_{[1:N]}^{\star} \leftarrow l_{\xi_{\mathrm{inv}}}\,\Delta \hat{\xi}_{[1:N]}^{j+1}$
    \State $\Delta \hat{u}_{[1:N]}^{\star} \leftarrow l_{u_{\mathrm{inv}}}\,\Delta \hat{u}_{[1:N]}^{j+1}$
    \State $\Delta \hat{s}^{\star} \leftarrow l_{s_{\mathrm{inv}}}\,\Delta \hat{s}^{j+1}$\vspace{1ex}
    \State $w_{[1:N-1]}^{\star} \leftarrow w_{[1:N-1]}^{j+1}$ \Comment{update dual variable}
    \end{algorithmic}
    \begin{flushleft}
        \textbf{Return:} $\Delta \hat{x}_{[1:N]}^{\star}$, $\Delta \hat{\xi}_{[1:N]}^{\star}$, $\Delta \hat{u}_{[1:N]}^{\star}$, $\Delta \hat{s}^{\star}$, $w^{\star}_{[1:N-1]}$
    \end{flushleft}
\end{algorithm}

\restoregeometry
    
\begin{algorithm}[H]
\small
\caption{Stopping Criterion Evaluation:\\[1ex]
\hphantom{\textbf{Algorithm 7}~}\scriptsize$\textsc{stopping}(\Delta\hat{x}_{[1:N]}^{j+1},\,\Delta\hat{\xi}_{[1:N]}^{j+1},\,\Delta\hat{u}_{[1:N]}^{j+1},\,\Delta\hat{s}^{j+1},\,w_{[1:N-1]}^{j+1}$,\\
\hphantom{\normalsize{\textbf{Algorithm 7}~}{\scriptsize{$\textsc{stopping}$\:}}}$\,\Delta\hat{x}_{[1:N]}^{j},\,\Delta\hat{\xi}_{[1:N]}^{j},\,\Delta\hat{u}_{[1:N]}^{j},\,\,\,\,\Delta\hat{s}^{j\hphantom{+1}}\!\!\!,\,w_{[1:N-1]}^{j},\,\epsilon_{\mathrm{abs}},\,\epsilon_{\mathrm{rel}})$}\label{alg:stopping}
    \vspace{0.25em}
    \begin{flushleft}
        \textbf{Inputs:}
        $\Delta\hat{x}_{[1:N]}^{j+1},\,\Delta\hat{\xi}_{[1:N]}^{j+1},\,\Delta\hat{u}_{[1:N]}^{j+1},\,\Delta\hat{s}^{j+1},\,w_{[1:N-1]}^{j+1}$,\\
\hphantom{\textbf{Inputs:}\ }$\,\Delta\hat{x}_{[1:N]}^{j},\,\Delta\hat{\xi}_{[1:N]}^{j},\,\Delta\hat{u}_{[1:N]}^{j},\,\,\,\,\Delta\hat{s}^{j\hphantom{+1}}\!\!\!,\,w_{[1:N-1]}^{j},\,\epsilon_{\mathrm{abs}},\,\epsilon_{\mathrm{rel}}$ 
    \end{flushleft}
    \begin{algorithmic}[1]
    \State $z^{j+1}_{\infty} \leftarrow \max\!\left\{\|\Delta\hat{x}_{[1:N]}^{j+1}\|_{\infty},\, \|\Delta\hat{\xi}_{[1:N]}^{j+1}\|_{\infty},\, \|\Delta\hat{u}_{[1:N]}^{j+1}\|_{\infty},\, |\Delta\hat{s}^{j+1}|\right\}$
    \State $z^{j\hphantom{+ 1}}_{\infty} \leftarrow \max\!\left\{\|\Delta\hat{x}_{[1:N]}^{j}\|_{\infty},\, \|\Delta\hat{\xi}_{[1:N]}^{j}\|_{\infty},\, \|\Delta\hat{u}_{[1:N]}^{j}\|_{\infty},\, |\Delta\hat{s}^{j}|\right\}$
    \State $z^{\Delta j}_{\infty} \,\,\leftarrow \max\!\left\{\|\Delta\hat{x}_{[1:N]}^{j+1} - \Delta\hat{x}_{[1:N]}^{j}\|_{\infty},\, \|\Delta\hat{\xi}_{[1:N]}^{j+1} - \Delta\hat{\xi}_{[1:N]}^{j}\|_{\infty},\, \|\Delta\hat{u}_{[1:N]}^{j+1} - \Delta\hat{u}_{[1:N]}^{j}\|_{\infty},\, |\Delta\hat{s}^{j+1} - \Delta\hat{s}^{j}|\right\}$\vspace{1ex}
    \State $r^{j+1}_{\infty} \leftarrow \|w_{[1:N-1]}^{j+1}\|_{\infty}$
    \State $r^{j\hphantom{+ 1}}_{\infty} \leftarrow \|w_{[1:N-1]}^{j}\|_{\infty}$
    \State $r^{\Delta j}_{\infty} \,\, \leftarrow \|w_{[1:N-1]}^{j+1} - w_{[1:N-1]}^{j}\|_{\infty}$\vspace{1ex}
    \If {$z^{\Delta j}_{\infty} \leq \epsilon_{\text{abs}} + \epsilon_{\text{rel}}\,\max\!\left\{z^{j+1}_{\infty},\,z^{j}_{\infty}\right\}$ \textbf{and} $r^{\Delta j}_{\infty} \leq \epsilon_{\text{abs}} + \epsilon_{\text{rel}}\,\max\!\left\{r^{j+1}_{\infty},\,r^{j}_{\infty}\right\}$}\vspace{1ex}
    \State $\textsc{terminate} \leftarrow \textsc{true}$\vspace{1ex}
    \Else
    \State $\textsc{terminate} \leftarrow \textsc{false}$\vspace{1ex}
    \EndIf\vspace{1ex}
    \end{algorithmic}
    \begin{flushleft}
        \textbf{Return:} \textsc{terminate}
    \end{flushleft}
    \vspace{0.25em}
\end{algorithm}
\section{Results}\label{sec:results}
\subsection{Offline Benchmarking}

We benchmark {\pipgc} against three state-of-the-art convex optimization solvers: {\ecos}, {\mosek}, and {\gurobi} \citep{domahidi2013ecos, mosek, gurobi}, by means of a lunar approach-phase test case, with a fixed final attitude quaternion, $q_{f}$. We use the absolute-variable version of the solver described in \citep{kamath2023customized}, which does not include row normalization. The solver parameter, $\lambda$, is manually tuned. The {\pipgc} solver is implemented via C code, generated using the \textsc{matlab} Coder \citep{matlab, coder}. The \textsc{yalmip} convex optimization modeling tool in \textsc{matlab} is used to parse the problem and interface with the off-the-shelf solvers \citep{lofberg2004yalmip}. All trials are run on a 2018 MacBook Pro with a 2.6 GHz 6-core Intel Core i7 processor and 16 GB of RAM.

\begin{figure}
    \centering\hspace{1.375em}
    \includegraphics[width=0.875\linewidth]{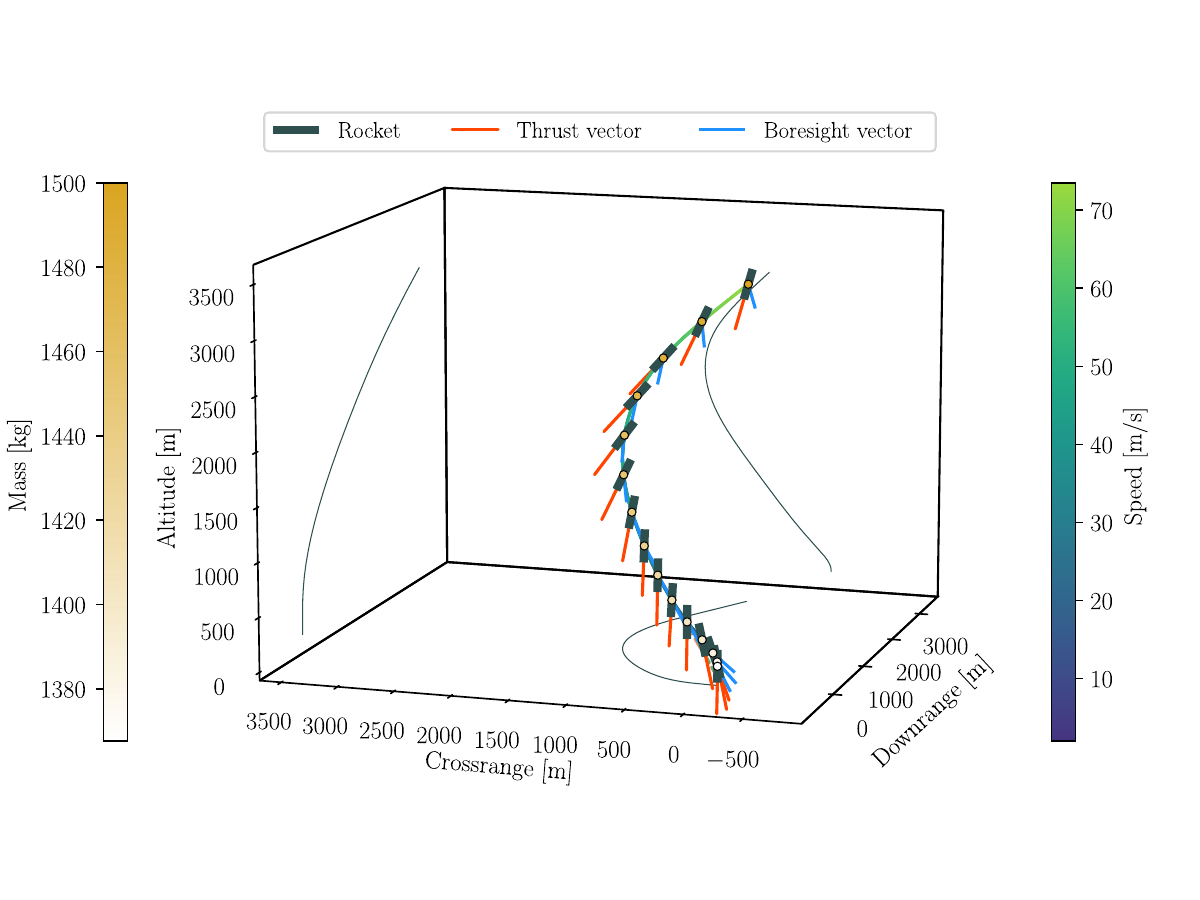}
    \vspace{-4em}
    \caption{The 3D landing trajectory obtained via SeCO in real-time ($N = 15$).}
    \label{fig:3D}
\end{figure}

\begin{figure}
    \centering
    \includegraphics[width=0.75\linewidth]{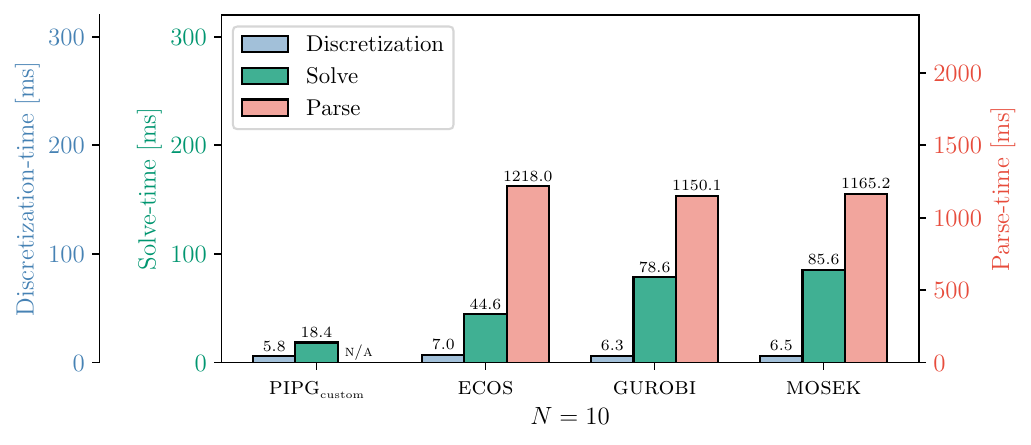}
    \includegraphics[width=0.75\linewidth]{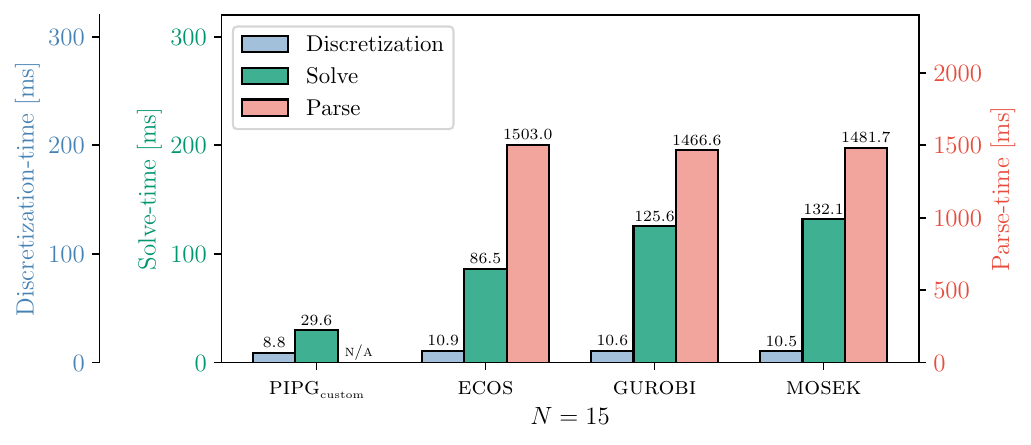}
    \includegraphics[width=0.75\linewidth]{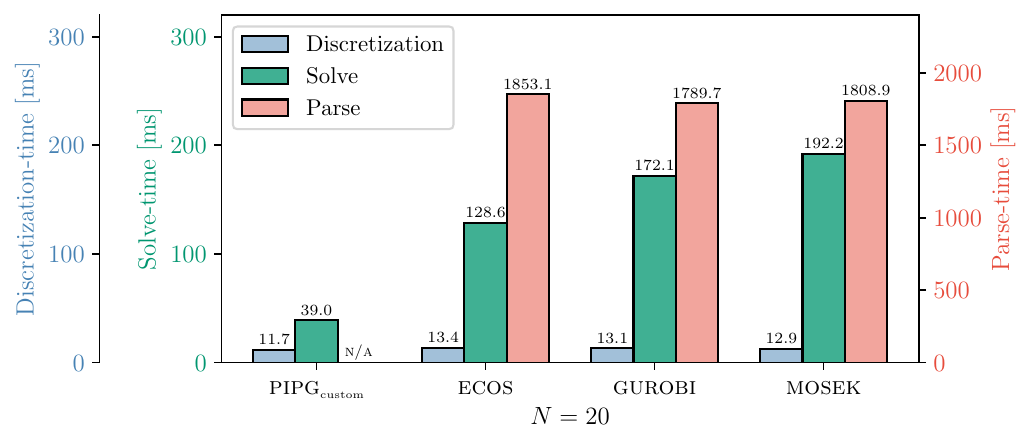}
    \includegraphics[width=0.75\linewidth]{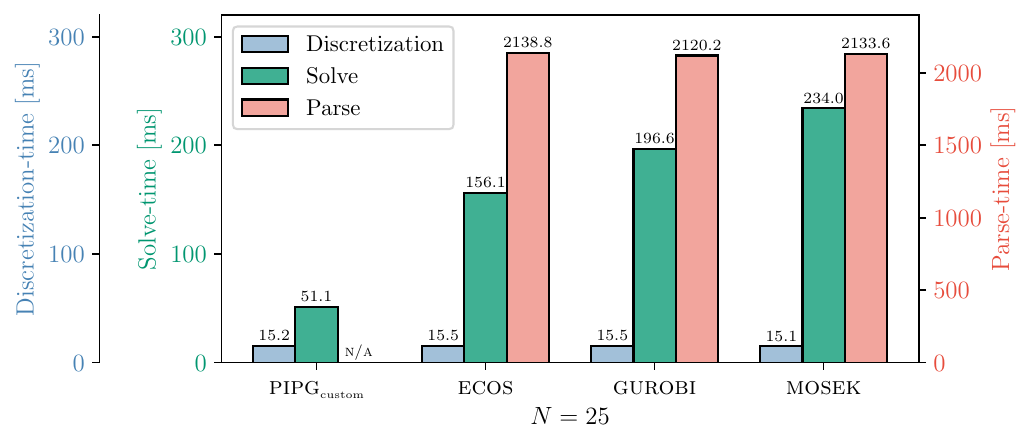}
    \vspace{2.5em}
    \caption{DQG benchmark test results.}
    \label{fig:benchmark}
\end{figure}

For consistency, the {\dqg} problem instance is set up such that each benchmarked solver solves the problem to a predetermined open-loop accuracy in exactly 5 {\seco} iterations. The entire {\dqg} problem is solved 100 times and the mean total (across all {\seco} iterations) discretization-, parse-, and solve-times are reported, as shown in Figures \ref{fig:benchmark}. The same procedure is carried out across 4 different problem sizes, representative of onboard guidance: $N \in \{10,\,15,\,20,\,25\}$, where $N$ is the number of discrete temporal nodes. The terminal position and velocity error tolerances (between the computed solution and the open-loop single-shot integrated trajectory) are set to $10$ m and $0.25$ m/s, respectively—similar to the tolerances chosen in \citep{reynolds2020dual}. Note that it is possible to significantly reduce parsing time for the other solvers for online execution \cite{reynolds2020real}, and the desktop parsing times are only reported for completeness; as such, the true performance comparison is between solve-times. The {\dqg} parameter values chosen for the benchmark test are given in Table \ref{tab:tab}. The 3-dimensional landing trajectory and the line-of-sight angle as a function of time (corresponding to $N = 15$), obtained via {\pipgc}, are shown in Figures \ref{fig:3D} and \ref{fig:los}, respectively.

We observe that the solution framework (\seco) itself leads to a speedup, regardless of the solver chosen, when compared with previously used SCP methods and solve-times reported in the literature \citep{reynolds2020real, reynolds2020dual, strohl2022implementation}. Further, {\pipg} is significantly faster than the solvers it is benchmarked against, as shown in Figure \ref{fig:timing}, and over an order of magnitude faster than the previously reported mean solve-time for {\dqg} \citep{reynolds2020dual} for the same problem size.

\begin{figure}
    \centering
    \includegraphics[width=0.65\linewidth]{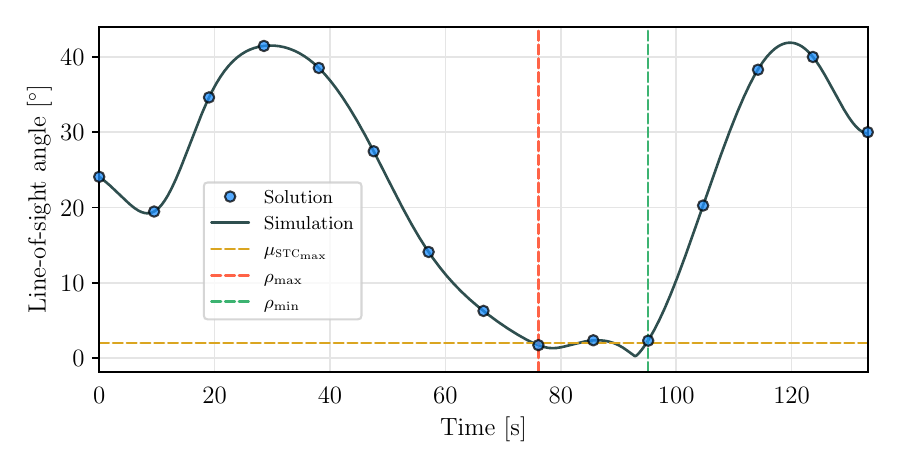}
    \vspace{-1em}
    \caption{The line-of-sight angle, which is constrained to be within 2$^\circ$ in the trigger window ($N = 15$).}
    \label{fig:los}
\end{figure}
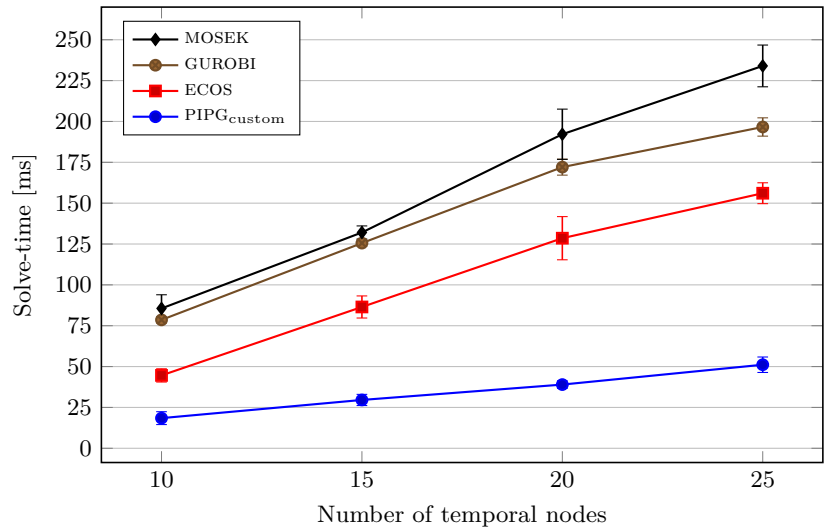
\begin{figure}
    \vspace{1.075em}
    \centering\hspace{-1.45em}
    \begin{tikzpicture}
    \begin{axis}[
    xlabel={Number of temporal nodes},
    ylabel={Solve-time [ms]},
    ymajorgrids=true,
    legend cell align={left},
    legend pos=north west,
    xtick=data,
    ytick distance={25},
    legend entries={\pipgc,\ecos,\gurobi,\mosek},
    width = 0.8*5.45in,
    height = 0.8*3.75in,
    label style={font=\small},
    tick label style={font=\small},
    legend style={font=\footnotesize},
    cycle list name=color,
    axis line style = {semithick},
    reverse legend,
    error bars/y dir=both, 
    error bars/y explicit  
    ]
    \addplot+[thick, error bars/.cd, error bar style={semithick}] table[x=x, y=y, y error=e] {data/PIPG_PGF_new.csv};
    \addplot+[thick, error bars/.cd, error bar style={semithick}] table[x=x, y=y, y error=e] {data/ECOS_PGF.csv};
    \addplot+[thick, error bars/.cd, error bar style={semithick}] table[x=x, y=y, y error=e] {data/Gurobi_PGF.csv};
    \addplot+[thick, mark=diamond*, error bars/.cd, error bar style={semithick}] table[x=x, y=y, y error=e] {data/Mosek_PGF.csv};
    \end{axis}
    \end{tikzpicture}
    \vspace{-0.25em}
    \caption{Solve-time comparison between the DQG-customized version of PIPG and three state-of-the-art convex optimization solvers. The error bars indicate three standard deviations ($\pm\!\>3\!\>\sigma$).}
    \label{fig:timing}
\end{figure}

\noindent\begin{table}[!htpb]

\centering

\begin{minipage}[t]{0.475\linewidth}
\centering

\setlength{\tabcolsep}{12.5pt}

{\renewcommand{\arraystretch}{1.25}\begin{tabular}{|ll|}
\hhline{|==|}
Parameter & Value \\ \hhline{|==|}
$g$ & $1.625$ m\,s$^{-2}$ \\ 
$g_{0}$ & $9.81$ m\,s$^{-2}$ \\ 
$I_{\mathrm{sp}_{\textsc{me}}}$ & $300$ s \\
$\alpha_{\textsc{me}}$ & $\frac{1}{I_{\mathrm{sp}_{\textsc{me}}} g_{0}}$ s\,m$^{-1}$\\
$I_{\mathrm{sp}_{\textsc{rcs}}}$ & $200$ s \\
$\alpha_{\textsc{rcs}}$ & $\frac{1}{I_{\mathrm{sp}_{\textsc{rcs}}} g_{0}}$ s\,m$^{-1}$\\
$\tau_{\max}$ & $50$ kg\,m$^{2}$\,s$^{-2}$ \\
$T_{\max}$ & $3000$ kg\,m\,s$^{-2}$ \\ 
$T_{\min}$ & $600$ kg\,m\,s$^{-2}$ \\
$\dot{T}_{\max}$ & $0.75\cdot(T_{\max} - T_{\min})$ kg\,m\,s$^{-3}$ \\
$\delta_{\max}$ & $5^{\circ}$ \\
$\dot{\delta}_{\max}$, $\dot{\phi}_{\max}$ & $5^{\circ}$\,s$^{-1}$ \\
$l_{\textsc{cm}}$ & $1$ m \\
$p_{\B}$ & $\left[0.5,\, 0,\, -\frac{\sqrt{3}}{2}\right]^{\top}$\\
$m_{i}$ & $1500$ kg \\
$m_{f}$ & $750$ kg \\
$J$ & $\operatorname{diag}\{4.2,\,4.2,\,0.6\}$ m$^{2}$ \\
\hhline{|==|}
\end{tabular}
}
\end{minipage}\hspace{0.275em}%
\begin{minipage}[t]{0.475\linewidth}

\setlength{\tabcolsep}{10pt}

{\renewcommand{\arraystretch}{1.25}\begin{tabular}{|ll|}
\hhline{|==|}
Parameter & Value \\ \hhline{|==|}
$\theta_{\max}$ & $90^{\circ}$ \\
$\omega_{\max}$ & $5^{\circ}$\,s$^{-1}$ \\
$v_{\max}$ & $90$ m\,s$^{-1}$ \vphantom{$\frac{1}{I_{\mathrm{sp}_{\textsc{me}}} g_{0}}$} \\
$h_{\min}$ & $100$ m \\
$\rho_{\max}$ & $1250$ m \\
$\rho_{\min}$ & $500$ m \\
$\theta_{\textsc{stc}_{\max}}$ & $20^{\circ}$ \vphantom{$\frac{1}{I_{\mathrm{sp}_{\textsc{rcs}}} g_{0}}$} \\
$\omega_{\textsc{stc}_{\max}}$ & $1^{\circ}$\,s$^{-1}$ \\
$v_{\textsc{stc}_{\max}}$ & $30$ m\,s$^{-1}$ \\
$\mu_{\textsc{stc}_{\max}}$ & $2^{\circ}$ \\
$r_{\I_{i}}$ & $\left[3000,\, 600,\, 3000\right]^{\top}$ m \\
$r_{\I_{f}}$ & $\left[0,\, 0,\, 100\right]^{\top}$ m \\
$v_{\I_{i}}$ & $\left[-60,\, 30,\, -30\right]^{\top}$ m\,s$^{-1}$ \\
$v_{z_{\I_{f}}}$ & $-2$ m\,s$^{-1}$ \\
$q_{i}$ & $\left[-0.15,\, 0.3,\, -1,\, 1\right]^{\top}$ (normalized) \vphantom{$\left.\frac{\sqrt{3}}{2}\right]$} \\
$q_{f}$ & $\left[0,\, 0,\, -1.25,\, 1\right]^{\top}$ (normalized) \\
$\omega_{\B_{i}}$ & $\left[0,\, 0,\, 0\right]^{\top}$ $^{\circ}$\,s$^{-1}$ \\[0.175em]
\hhline{|==|}
\end{tabular}
}
\end{minipage}
\vspace{0.75em}
\caption{The DQG parameter values chosen for the solver benchmark test.}
\label{tab:tab}

\end{table}

\newpage
\subsection{Onboard (Hardware-in-the-Loop) Testing}

\begin{figure}[H]
    \centering\hspace{1.375em}
    \includegraphics[width=0.875\linewidth]{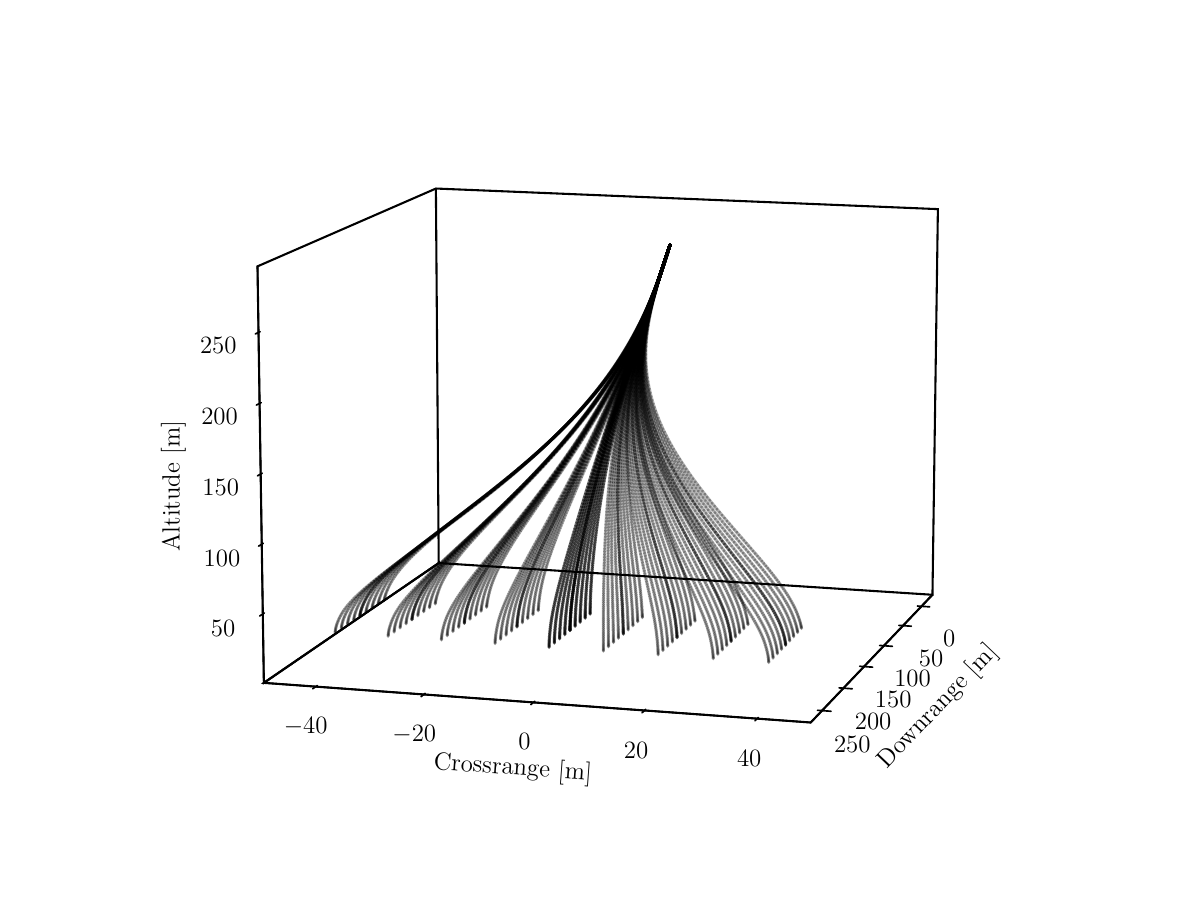}
    \vspace{-3em}
    \caption{Hazard-avoidance divert trajectories computed onboard the NASA SPLICE Descent and Landing Computer (DLC) in a hardware-in-the-loop setting.}
    \label{fig:hitl_mc}
\end{figure}

We consider a terrestrial rocket landing mission scenario for an upcoming closed-loop (\dqg{}-in-the-loop) rocket landing flight test campaign \citep{Mendeck_SPLICE_2023, Mendeck_SPLICE_2024}, and solve the problem, for $100$ divert sites on a uniform grid, as shown in Figure \ref{fig:hitl_mc}. The custom solver used for this application, based on deviation variables, is detailed in Section \ref{sec:custom}. We perform row-normalization, and manually tune the solver parameter, $\lambda$. The problem formulation is identical to Problem \ref{subsec:cont_time_ocp}, with a few modifications, such as an independent thrust and torque model (without gimbaling of the rocket engine in guidance and with control allocation handled outside of guidance), the inclusion of aerodynamic forces, independent component-wise torque bounds, a state-triggered glideslope constraint to replace the state-triggered tilt constraint, and imposition of the initial condition constraint on the true state (as opposed to the virtual state). See \citep{doll2025hardware} for more details on the problem formulation. The {\pipgc} solver is implemented via C code, which, again, is generated using the \textsc{matlab} Coder \citep{matlab, coder}. The solver is executed onboard the NASA SPLICE Descent and Landing Computer (DLC), which consists of a cluster of $4$ ARM Cortex A53 processors, on which the flight software runs \citep{rutishauser2023system}.

For the terrestrial landing scenario considered in \cite{fritz2022post} for the suborbital flight tests of the Blue Origin New Shepard reusable launch vehicle (with \dqg{} executed in an open-loop, onboard the DLC), with the older SCP algorithm and a customized version of the IPM-based subproblem solver, \bsocp{}, the $10$-node version took an average (over $3$ runs) of $2.65$ seconds, with all $3$ runs taking $4$ SCP iterations each.

For the lunar landing scenario considered in \citep{strohl2022implementation}, with the older SCP algorithm and the IPM-based subproblem solver, \bsocp{} \citep{dueri2014automated, dueri2017customized}, the $10$-node version took an average (over $100$ runs) of $5.85$ seconds, with $97$ runs taking $4$ SCP iterations each, and the remaining $3$ runs taking $5$ SCP iterations each. The $20$-node version, on the other hand, took an average (again, over $100$ runs) of $11.36$ seconds, with $27$ runs taking $4$ SCP iterations each, and the remaining $73$ runs taking $3$ SCP iterations each. Neither of these versions met the SPLICE goal of a guidance update-rate of $1$ second, or even the SPLICE requirement of a guidance update-rate of $3$ seconds.

In contrast, the custom solver proposed in this work, applied to the terrestrial landing scenario considered in \cite{doll2025hardware}, with $15$ nodes, took an average (over $100$ divert scenarios) of $0.5887$ seconds, with all $100$ runs taking $5$ SCP (\seco{}) iterations each. Note that this not only meets, but exceeds both the SPLICE requirement of a guidance update-rate of $3$ seconds and the SPLICE goal of a guidance update-rate of $1$ second, by a significant margin. Further, we note that this marks the first time in the duration of the NASA SPLICE program that the $1$-second goal has been achieved onboard the DLC.

\begin{figure}[H]
    \centering
    \includegraphics[width=0.925\linewidth]{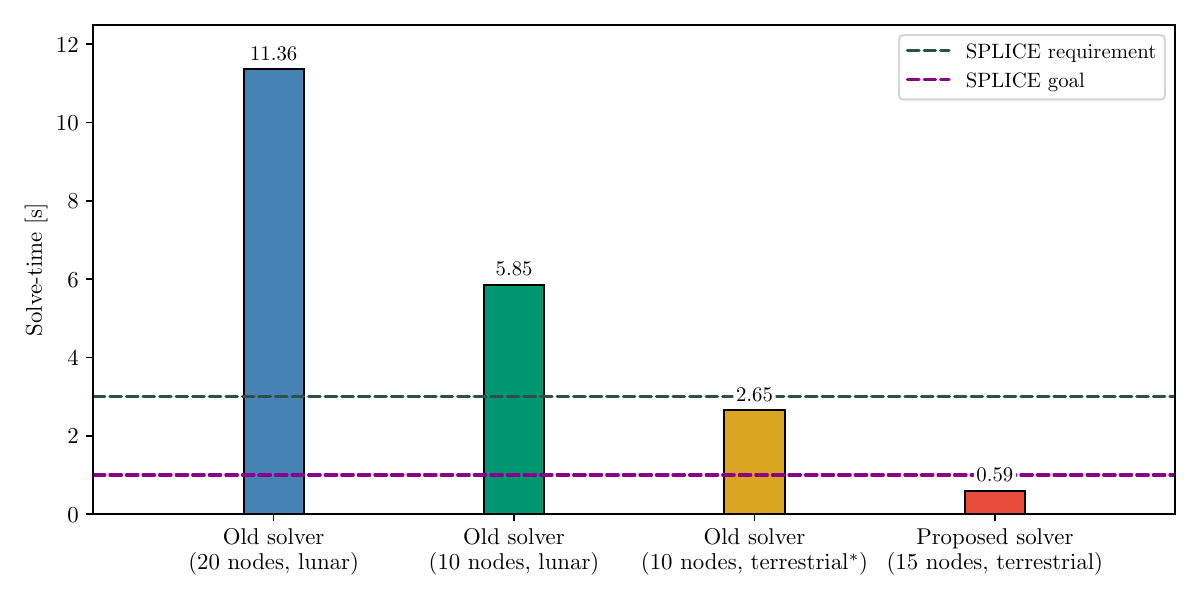}
    \caption{Average solve-times onboard the NASA SPLICE Descent and Landing Computer (DLC). The ``old solver'' refers to the previously-used SCP algorithm \cite{reynolds2020dual} with the BSOCP convex subproblem solver \cite{dueri2017customized}. From the left, the first two bars (averaged over 100 runs) are from \cite{strohl2022implementation} with generic BSOCP, the third bar (averaged over 3 runs) is from \cite{fritz2022post} with customized BSOCP ($^{*}$in-flight), and the rightmost bar (averaged over 100 runs) corresponds to the proposed solver with customized PIPG (this work), which meets both the SPLICE requirement and the SPLICE goal for the guidance update-rate, i.e., solve-time.}
    \label{fig:onboard}
\end{figure}

These results are presented in Figure \ref{fig:onboard}. We note that this does not represent a direct comparison, owing to the differing problem formulations and parameters considered in the preceding tests between the different solvers. That said, given the similarity in mission complexity and the fact that all solvers were executed on the same computing platform (the DLC), we conclude that our proposed solver is roughly $5$ to $10$ times faster than the old solver.

\vspace{-2.25em}

\vspace{3em}

\newpage
\section{Conclusions}

Sequential conic optimization ({\seco}) combines sequential convex programming (SCP) with first-order conic optimization to solve difficult trajectory optimization problems, such as the dual quaternion-based 6-DoF powered-descent guidance ({\dqg}) problem, in real-time. First-order optimization solvers, such as {\pipg}, are attractive for: (i) real-time applications (given their execution speed); (ii) implementation onboard resource-constrained systems (owing to the small footprint of the resulting codebase); and, (iii) verification and validation (due to their reliance on simple computations). Recent advances have enabled this class of algorithms to match (and even out-perform) solvers based on interior-point methods (IPMs). Further, {\pipg} is amenable to warm-starting and performance-efficient customization for trajectory optimization problems.

We formulate the nonconvex {\dqg} problem---with mission-critical constraints---in compliance with the {\seco} framework, and solve it using $\pipgc$, a custom first-order conic optimization solver developed for this application, in conjunction with a customized preconditioning algorithm. This solver is able to solve the entire nonconvex problem in a matter of milliseconds, and is much faster than other state-of-the-art convex optimization solvers across varying problem sizes.

Finally, we demonstrate, by means of hardware-in-the-loop testing onboard the NASA SPLICE Descent and Landing Computer (DLC), that the resulting algorithm can generate trajectories fast enough in terms of satisfying NASA's guidance update-rate requirements for hazard detection and avoidance (HDA) maneuvers for autonomous precision rocket-landing.
\section*{Acknowledgements}

The authors thank the members of the Autonomous Controls Laboratory (ACL) at the University of Washington, especially Dayou Luo and Samet Uzun, for the discussions on solver development and acceleration, Govind Chari, for a detailed review of the manuscript, and Benjamin Chung, for the discussions on sparse linear algebra. We also thank the members of the Flight Mechanics and Trajectory Design branch (EG5) at the NASA Johnson Space Center, especially Breanna Johnson, Dan Matz, and Ron Sostaric, for their valuable guidance, insight, and many helpful discussions. The authors give their special thanks to the co-developers of the original {\dqg} algorithm, Miki Szmuk and Danylo Malyuta, for their ongoing support. This research was supported by NASA grant NNX17AH02A and was partially carried out at the NASA Johnson Space Center; Government sponsorship is acknowledged.

\newpage

\bibliography{references}

\newpage
\section*{Appendix}
\appendix
\label{appendix}
\subsection{Quaternion Algebra}
\subsubsection{Unit Quaternions}
\[\q{a} = \left[\q{a}_{v}^{\top},\,a_{4}\right]^{\top},\ \q{b} = \left[\q{b}_{v}^{\top},\,b_{4}\right]^{\top} \in\,\mathbb{R}_{u}^{4} \defeq \left\{\q{q} \infour\ \Big\vert\ \q{q}^{\top}\q{q} = 1\right\}\]
where \[\q{a}_{v} = \left[a_{1},\,a_{2},\,a_{3}\right]^{\top},\ \q{b}_{v} = \left[b_{1},\,b_{2},\,b_{3}\right]^{\top} \inthree\ \text{and}\ a_{4}, b_{4} \inone\]
\underline{Note}: All quaternions are in accordance with the scalar-last convention.
\subsubsection{Conjugation}
\[\q{a}^{*} \defeq \left[-\q{a}_{v}^{\top},\,a_{4}\right]^{\top}\]
\subsubsection{Skew-Symmetric Matrix Operator}
\[\q{a}_{v}^{\times} \defeq \begin{pmatrix}
       0 & -a_{3} & a_{2}\\
       a_{3} & 0 & -a_{1}\\
       -a_{2} & a_{1} & 0
    \end{pmatrix}\]
\subsubsection{\texorpdfstring{$SO(4)$}{SO4} Matrix Operators}
\[[\q{a}]_{\otimes}\defeq\begin{pmatrix}
a_{4} I_{3}+\q{a}_{v}^{\times} & \q{a}_{v} \\
-\q{a}_{v}^{\top} & a_{4}
\end{pmatrix}\]
\[[\q{b}]_{\otimes}^{*}\defeq\begin{pmatrix}
b_{4} I_{3}-\q{b}_{v}^{\times} & \q{b}_{v} \\
-\q{b}_{v}^{\top} & b_{4}
\end{pmatrix}\]
\subsubsection{Multiplication}
\vspace{-1em}
\begin{align*}
\q{a}\otimes\q{b}\ &{\defeq} \left[a_{4} \q{b}_{v} + b_{4} \q{a}_{v} + \q{a}_{v}\times\q{b}_{v},\,a_{4} b_{4} - \q{a}_{v}^{\top}\q{b}_{v}\right]^{\top}\\
&= [\q{a}]_{\otimes}\q{b}\\
&= [\q{b}]_{\otimes}^{*}\q{a}
\end{align*}
\subsubsection{Cross Product}
\[\q{a}\oslash\q{b} \defeq \left[a_{4} \q{b}_{v} + b_{4} \q{a}_{v} + \q{a}_{v}\times\q{b}_{v},\,0\right]^{\top}\]
\subsection{Dual Quaternion Algebra}
\subsubsection{Unit Dual Quaternions}
\[\dq{a} = \left[\q{a}_{1}^{\top},\,\q{a}_{2}^{\top}\right]^{\top},\ \dq{b} = \left[\q{b}_{1}^{\top},\,\q{b}_{2}^{\top}\right]^{\top} \in\,\mathbb{R}_{u}^{8} \defeq \left\{\dq{q} = \left[\q{q}_{1}^{\top},\,\q{q}_{2}^{\top}\right]^{\top} \ineight\ \Big\vert\ \q{q}_{1}^{\top}\q{q}_{1} = 1~\text{and}~\q{q}_{1}^{\top}\q{q}_{2} = 0,\ \q{q}_{1}, \q{q}_{2} \infour\right\}\]
\subsubsection{Conjugation}
\[\dq{a}^{*} \defeq \begin{pmatrix}\q{a}_{1}^{*}\\\q{a}_{2}^{*}\end{pmatrix} \label{eq:dq_conj}\]
\subsubsection{Multiplication Matrix Operators}
\[[\dq{a}]_{\otimes} \defeq \begin{pmatrix}
[\q{a}_{1}]_{\otimes} & 0_{4\times4} \\
[\q{a}_{2}]_{\otimes} & [\q{a}_{1}]_{\otimes}
\end{pmatrix}\]
\[[\dq{b}]_{\otimes}^{*} \defeq \begin{pmatrix}
[\q{b}_{1}]_{\otimes}^{*} & 0_{4\times4} \\
[\q{b}_{2}]_{\otimes}^{*} & [\q{b}_{1}]_{\otimes}^{*}
\end{pmatrix}\]
\subsubsection{Multiplication}
\[\dq{a}\otimes\dq{b} \defeq [\dq{a}]_{\otimes}\dq{b} = [\dq{b}]_{\otimes}^{*}\dq{a} \label{eq:dq_mul}\]
\subsubsection{Cross Product}
\[\dq{a}\oslash\dq{b} \defeq \left[\q{a}_{1}\oslash\q{b}_{1},\, \q{a}_{1}\oslash\q{b}_{2} + \q{a}_{2}\oslash\q{b}_{1}\right]^{\top}\]

\end{document}